\documentclass[a4paper,12pt]{article}
\usepackage[T1]{fontenc}
\usepackage[utf8]{inputenc}
\usepackage{theorem} 
\usepackage{amsmath} 
\usepackage{amsfonts}
\usepackage{amscd} 
\usepackage{amssymb} 
\usepackage{enumerate}
\usepackage{mathrsfs}
\usepackage{url}
\usepackage[francais]{babel}
\theoremstyle{break}
\newtheorem{thm}{Th\'eor\`eme}[section]
\newtheorem{prop}[thm]{Proposition}
\newtheorem{lemme}[thm]{Lemme}
\newtheorem{cor}[thm]{Corollaire}
\newtheorem{question}[thm]{Question}
\newtheorem{defi}[thm]{D\'efinition}
\newtheorem{notas}[thm]{Notations}
\newtheorem{nota}[thm]{Notation}
\newenvironment{demo}
{\par {\em D\'emonstration :}{\hspace{0.02\textwidth}}}{\hfill $\Box$\par}

\newenvironment{rem}
{\refstepcounter{thm}\noindent{\em Remarque \arabic{section}.\arabic{thm} :}{\hspace{0.02\textwidth}}}{\hfill $\Box$\par}

\makeatletter
\@addtoreset{equation}{section}

\makeatother
\newcommand{\termin}[1]{{\em #1}}
\newcommand{\C}{\mathbf{C}}
\newcommand{\G}{\mathbf{G}}
\newcommand{\N}{\mathbf{N}}
\newcommand{\R}{\mathbf{R}}
\newcommand{\Z}{\mathbf{Z}}
\newcommand{\bT}{\boldsymbol{T}}
\newcommand{\ba}{\boldsymbol{a}}
\newcommand{\bd}{\boldsymbol{d}}
\newcommand{\be}{\boldsymbol{e}}
\newcommand{\bbf}{\boldsymbol{f}}
\newcommand{\bm}{\boldsymbol{m}}
\newcommand{\bn}{\boldsymbol{n}}
\newcommand{\bx}{\boldsymbol{x}}
\newcommand{\by}{\boldsymbol{y}}
\newcommand{\bz}{\boldsymbol{z}}
\newcommand{\bnu}{\boldsymbol{\nu}}
\newcommand{\cC}{{\cal C}}
\newcommand{\cD}{{\cal D}}
\newcommand{\cE}{{\cal E}}
\newcommand{\cF}{{\cal F}}
\newcommand{\cT}{{\cal T}}
\DeclareMathAlphabet{\eulercal}{U}{eus}{m}{n}
\newcommand{\ecA}{{\eulercal A}}
\newcommand{\ecD}{{\eulercal D}}
\newcommand{\ecE}{{\eulercal E}}
\newcommand{\ecF}{{\eulercal F}}
\newcommand{\ecG}{{\eulercal G}}
\newcommand{\ecH}{{\eulercal H}}
\newcommand{\ecI}{{\eulercal I}}
\newcommand{\ecL}{{\eulercal L}}
\newcommand{\ecN}{{\eulercal N}}
\newcommand{\ecO}{{\eulercal O}}
\newcommand{\ecT}{{\eulercal T}}
\DeclareMathAlphabet{\beulercal}{U}{eus}{b}{n}
\newcommand{\becD}{{\beulercal D}}
\newcommand{\becE}{{\beulercal E}}
\newcommand{\becG}{{\beulercal G}}
\newcommand{\frD}{{\mathfrak D}}
\newcommand{\bfrD}{\boldsymbol{{\mathfrak D}}}
\newcommand{\frE}{{\mathfrak E}}
\newcommand{\bfrE}{\boldsymbol{{\mathfrak E}}}
\newcommand{\frF}{{\mathfrak F}}
\newcommand{\longeq}{\Longleftrightarrow}
\newcommand{\longto}{\longrightarrow}
\newcommand{\isom}{\overset{\sim}{\to}}
\newcommand{\longisom}{\overset{\sim}{\longto}}
\DeclareMathOperator{\Cox}{Cox}
\DeclareMathOperator{\Pic}{Pic}
\DeclareMathOperator{\Div}{Div}
\DeclareMathOperator{\rg}{rg}
\DeclareMathOperator{\Ker}{Ker}
\DeclareMathOperator{\im}{Im}
\DeclareMathOperator{\Hom}{Hom}
\DeclareMathOperator{\Min}{Min}
\DeclareMathOperator{\Max}{Max}
\DeclareMathOperator{\Spec}{Spec}
\DeclareMathOperator{\pgcd}{pgcd}
\let\leq\leqslant
\let\geq\geqslant
\newcommand{\sumu}[1]{\underset{#1}{\sum}}
\newcommand{\produ}[1]{\underset{#1}{\prod}}
\newcommand{\cupu}[1]{\underset{#1}{\cup}}
\newcommand{\oplusu}[1]{\underset{#1}{\oplus}}
\newcommand{\otimesu}[1]{\underset{#1}{\otimes}}
\newcommand{\disju}[1]{\underset{#1}{\bigsqcup}}
\newcommand{\Maxu}[1]{\underset{#1}{\Max}}
\newcommand{\Minu}[1]{\underset{#1}{\Min}}
\newcommand{\simu}[1]{\underset{#1}{\sim}}
\newcommand{\eps}{\varepsilon}
\newcommand{\vide}{\varnothing}
\newcommand{\eqdef}{\overset{\text{{\tiny{d\'ef}}}}{=}}
\newcommand{\ind}{\mathbf{1}}
\newcommand{\wt}{\widetilde}
\newcommand{\acc}[2]{\left\langle #1\, ,\,#2 \right\rangle} 
\newcommand{\norm}[1]{\left|\left| #1 \right|\right|} 
\newcommand{\abs}[1]{\left| #1 \right|} 
\newcommand{\card}[1]{\# #1} 
\newcommand{\ecOu}[1]{\underset{#1}{\ecO}}
\renewcommand{\div}{\text{div}}
\newcommand{\classe}[1]{\left[#1\right]} 
\newcommand{\idex}{\mathscr{I}_X} 
\newcommand{\idexh}{\mathscr{I}_{X,\text{homog}}} 
\newcommand{\scan}[1]{s_{#1}} 
\newcommand{\inv}{\times}
\newcommand{\courbe}{\mathscr{C}}
\newcommand{\cDi}{\cD_i}
\newcommand{\sep}{\,\text{\textnormal{sep}}}
\newcommand{\ceff}{C_{\text{eff}}}
\newcommand{\TNS}{T_{\text{NS}}}
\newcommand{\tors}{\cT}
\newcommand{\Homog}{\ecH}
\newcommand{\Homogs} {\ecH^{\bullet}}
\newcommand{\diveffc} {\Div_{\text{eff}}(\courbe)}
\newcommand{\Zp} {Z_{0}}
\newcommand{\Zpr}[1]{Z_{0,#1}}
\newcommand{\wtZpr}[1]{\wt{Z}_{0,#1}}
\DeclareMathOperator{\fact}{fact}
\DeclareMathOperator{\dens}{dens}
\DeclareMathOperator{\ddiv}{div}
\newcommand{\tpiczc} {\wt{\Pic}^0(\courbe)}
\newcommand{\piczcixe}{\wt{\Pic}^0(\courbe)^I_{X,\becE}}
\newcommand{\piczcSe} {\wt{\Pic}^0(\courbe)^7_{S,\becE}}
\newcommand{\nix}{\N^{\,I}_X}
\title{Comptage de courbes sur le plan projectif \'eclat\'e en trois
  points align\'es}
\author{David Bourqui}
\date{}
\begin{document}
\maketitle
{\small{ \textbf{R\'esum\'e} : Nous \'etablissons une version de la
  conjecture de Manin pour le plan projectif \'eclat\'e en trois
  points align\'es, le corps de base \'etant un corps global de
  caract\'eristique positive.}}
\par
{\small{ \textbf{Abstract} : We prove a version of Manin's conjecture 
for the projective plane blown up in three collinear points, the base
field being a global field of positive characteristic.}}
\vskip1cm

Vers la fin des ann\'ees 1980, une s\'erie de questions a \'et\'e soulev\'ee par Manin et ses collaborateurs 
(cf. \cite{BaMa:pts_rat,FMT}) sur le comportement
asymptotique du nombre de points de hauteur born\'ee des vari\'et\'es de
Fano d\'efinies sur un corps global, ce qui a initi\'e de nombreux
travaux. Ces derni\`eres ann\'ees, ces questions ont notamment \'et\'e \'etudi\'ees de mani\`ere
intensive pour les surfaces de del Pezzo g\'en\'eralis\'ees
d\'efinies sur un corps de nombres, en 
suivant la strat\'egie initi\'ee par Salberger dans \cite{Sal:tammes},
dont le principe consiste \`a relever le d\'ecompte des points de
hauteur born\'ee \`a un certain torseur sous un tore au-dessus
de la vari\'et\'e \'etudi\'ee (en g\'en\'eral un torseur universel).
Dans de nombreux cas de surfaces, le probl\`eme de comptage \og relev\'e\fg\ peut alors
se traiter par des techniques de th\'eorie analytique des nombres.
Un des points cruciaux pour l'efficacit\'e de la m\'ethode est que l'on
dispose d'\'equations explicites pour les torseurs
consid\'er\'es. 
Voici un tr\`es bref aper\c cu des r\'esultats obtenus dans ce contexte (nous
renvoyons au survol \cite{Br:manin:conjecture:dim:2} pour plus de
d\'etails). Le premier exemple de surface de del Pezzo non torique \`a
avoir \'et\'e trait\'e est d\^u \`a de la Bret\`eche  : il s'agit de l'\'eclat\'e du plan projectif en 4
points en position g\'en\'erale, {\it i.e.} de la surface de del Pezzo (lisse)
de degr\'e 5 (\cite{dlB:duke}). Par la suite, un certain nombre d'exemples de surfaces de del Pezzo
singuli\`eres de degr\'e 3, 4, 5 ou 6 ont \'et\'e trait\'es, essentiellement par de la Bret\`eche,
Browning et Derenthal ({\it cf.} par exemple
\cite{dlBBD,dlBBr:sdP4I}). L'obtention d'un \'equivalent asymptotique
pour le nombre de points de hauteur born\'ee dans le 
cas des surfaces cubiques lisses semble notamment hors de port\'ee pour l'instant.

Dans ce texte, nous \'etudions le probl\`eme dans le cas o\`u la
surface $X$
consid\'er\'ee  est le plan projectif \'eclat\'e en trois points align\'e, le corps de base \'etant
un corps global de caract\'eristique positive. Nous \'etablissons une
version de la conjecture de Manin dans ce cas (th\'eor\`eme \ref{thm:resprinc}).
Dans ce cadre, le probl\`eme a une reformulation g\'eom\'etrique simple : 
on cherche \`a \'evaluer le nombre de morphismes d'une courbe fix\'ee
$\courbe$ vers $X$ de degr\'e donn\'e, quand le degr\'e devient grand.
Soulignons que la courbe $\courbe$ est suppos\'ee de genre quelconque
dans notre \'etude, alors que les travaux cit\'es pr\'ec\'edemment sur les surfaces de del
Pezzo g\'en\'eralis\'ees se limitent en g\'en\'eral au
cas o\`u le corps de base est le corps des rationnels. 

Nous suivons \'egalement la strat\'egie de Salberger : la premi\`ere \'etape de la d\'emonstration
consiste \`a relever le probl\`eme de comptage au torseur universel au-dessus de $X$.
Nous d\'ecrivons cette \'etape dans la section \ref{sec:relev},
o\`u nous nous pla\c cons en fait dans le cadre plus g\'en\'eral d'une vari\'et\'e $X$ dont l'anneau
de Cox est suppos\'e de type fini.
De fa\c con informelle, ceci permet de ramener le d\'enombrement des morphismes de
$\courbe$ vers $X$ de degr\'e donn\'e \`a celui de certaines familles
de sections globales de fibr\'es en droites, sections
astreintes \`a satisfaire deux types de conditions : elles doivent
satisfaire les \og m\^emes\fg\
\'equations que le torseur universel et leurs diviseurs doivent satisfaire certaines conditions de
coprimalit\'e (conditions enti\`erement explicites en termes des donn\'ees d\'ecrivant
le torseur universel). On se reportera \`a la proposition
\ref{prop:descr:mor} pour un \'enonc\'e pr\'ecis.
Une technique standard d'inversion de M\"obius 
permet d'\og oublier\fg\ la condition de coprimalit\'e. Le plus ardu est de
tenir compte des \'equations satisfaites par les sections dans le d\'enombrement. 
\`A cet \'egard, le cas le plus simple appara\^\i t comme \'etant celui des vari\'et\'es toriques, o\`u il n'y a
pas d'\'equations. On est alors ramen\'e \`a estimer la dimension
de certains espaces de sections globales, ce qui peut se faire via le
th\'eor\`eme de Riemann-Roch, et m\`ene \`a la d\'emonstration de la conjecture
de Manin dans ce cas. Ceci est fait dans \cite{Bou:vtetor}.
Dans cet article nous expliquons comment traiter la premi\`ere condition
dans le cas o\`u $X$ est le plan projectif \'eclat\'e en trois points align\'es, o\`u il n'y a qu'une
\'equation pour le torseur universel, qui plus est particuli\`erement simple. Les techniques utilis\'ees sont
\'el\'ementaires et susceptibles de s'adapter \`a d'autres surfaces de del Pezzo
g\'en\'eralis\'ees dont le torseur universel est donn\'e par une seule
\'equation ({\it cf.} \cite{Der:sdp:ut:hyp} pour la classification de ces
surfaces).

Nous terminons cette introduction par quelques remarques : 
dans le cas o\`u le corps de base est le corps des rationnels,
l'\'etude du nombre de points de hauteur born\'ee sur le plan projectif
\'eclat\'e en trois points align\'es est trait\'e via l'usage du torseur universel par
Browning dans \cite{Br:manin:conjecture:dim:2}. Dans ce cas, le
r\'esultat d\'ecoule aussi d'un th\'eor\`eme plus g\'en\'eral de Chambert-Loir et
Tschinkel sur la validit\'e des conjectures de Manin pour les
compactifications \'equivariantes d'espaces affines d\'efinies sur un
corps de nombres ({\it cf.} \cite{CLT_vect1,CLT_vect3}). Ces derniers auteurs
utilisent des techniques d'analyse harmonique. Il est probable que leurs
arguments s'adaptent dans le cas d'un corps global de caract\'eristique
positive, mais ceci resterait \`a mettre en \oe uvre.

\section{Rel\`evement du probl\`eme de comptage au torseur universel}\label{sec:relev}
\subsection{Quelques rappels sur la th\'eorie des anneaux de Cox}\label{subsec:rappels:cox}
Nous faisons quelques rappels sur la th\'eorie des anneaux de Cox,
initi\'ee par Cox dans \cite{Cox:hom_coo_ring} dans le cas des vari\'et\'es toriques.
On peut la voir comme une g\'en\'eralisation des coordonn\'ees homog\`enes sur
les espaces projectifs. Hassett et Tschinkel ont montr\'e qu'elle
fournissait un outil efficace pour la d\'etermination explicite des \'equations 
de certains torseurs universels ({\it cf.} \cite{Has:eq:ut:cox:rings}).

Soit $k$ un corps et $X$ une vari\'et\'e projective, lisse 
et g\'eom\'etriquement int\`egre d\'efinie sur $k$.
On suppose que le groupe de Picard de $X$ est libre de rang fini et
d\'eploy\'e, {\it i.e.} $\Pic(X_{\bar{k}})$ co\"\i ncide avec
$\Pic(X_{k^{\sep}})$
et l'action du groupe de Galois absolu est triviale.
On suppose en outre que l'anneau de Cox de $X$ ({\it cf.} \cite{Has:eq:ut:cox:rings}), not\'e $\Cox(X)$,
est de type fini. 
Soit $(s_i)_{i\in I}$ une famille finie de sections globales (non constantes) 
qui engendrent $\Cox(X)$.
Pour $i\in I$, soit $\cDi$ le diviseur des z\'eros de $s_i$. 
On note $D_i=\classe{\cDi}$ sa classe dans $\Pic(X)$.
 On d\'efinit une $\Pic(X)$-graduation
sur $k[(s_i)_{i\in I}]$ en posant $\deg(s_i)=D_i$.

Soit 
$\idex$ l'id\'eal 
$\Pic(X)$-homog\`ene noyau du morphisme naturel $k[s_i]\to \Cox(X)$,
de sorte qu'on a un isomorphisme
$
\Cox(X)\isom k[(s_i)]/\idex
$.

Soit $X_0$ le compl\'ementaire de la r\'eunion des $\cD_i$ pour $i\in I$.
On a la suite exacte classique
\begin{equation}\label{eq:suite_exacte_pic_tor}
0\longto k[X_0]^{\times}/k^{\times}\longto 
\oplusu{i\in I} \Z\,\cDi\longto \Pic(X) \longto 0.
\end{equation}
D'apr\`es le lemme de Rosenlicht, $k[X_0]^{\inv}/k^{\inv}$
est un $\Z$-module libre de rang fini.
\begin{nota}
On note $N_X$ le $\Z$-module $k[X_0]^{\inv}/k^{\inv}$.
\end{nota}
Dans la suite, 
on choisit arbitrairement une identification de $N_X$ \`a un sous-groupe de $k[X_0]^{\inv}$.

\begin{rem}\label{rem:cone:effectif}
Comme les sections $s_i$ engendrent $\Cox(X)$, 
le c\^one effectif de $X$, not\'e $\ceff(X)$,
 est l'image du c\^one $\oplusu{i\in I} \R_{\geq 0}\,\cDi$.
En particulier, il est de type fini.
\end{rem}

Le tore de N\'eron-Severi $\TNS\eqdef\Hom(\Pic(X),\G_m)$ agit
naturellement sur $k[s_i]$ via 
\begin{equation}
t.s_i=t(D_i).s_i
\end{equation}
et cette action induit une action de $\TNS$ sur $\Cox(X)$.

Soit $D$ une classe ample de $\Pic(X)$ et  
$\tors_X$ l'ensemble des points
de $\Spec(\Cox(X))$ semi-stables vis-\`a-vis de la $\TNS$-lin\'earisation sur le fibr\'e trivial de
$\Spec(\Cox(X))$ donn\'ee par $D$ (vu comme caract\`ere de $\TNS$). 
C'est un ouvert non vide $\TNS$-stable.
La th\'eorie g\'eom\'etrique des invariants permet alors de montrer
que le quotient g\'eom\'etrique de $\tors_X$ par $\TNS$ existe et s'identifie
naturellement \`a $X$. On montre en outre que $\tors_X\to X$ repr\'esente
l'unique classe de torseurs universels au-dessus de $X$ ({\it cf.} \cite{Has:eq:ut:cox:rings,hukeel:mori}).

Ainsi, si $D$ est une classe ample de $\Pic(X)$,
$\tors_X$ est le compl\'ementaire du ferm\'e dont l'id\'eal est engendr\'e par les 
\'el\'ements de $\Cox(X)$ homog\`enes de degr\'es $(m\,D)_{m\geq 1}$.
En consid\'erant des g\'en\'erateurs monomiaux de cet id\'eal,
on voit que, pour un certain ensemble $\ecI_X$ de parties de $I$, 
$\tors_X$ est la r\'eunion des ouverts de $\Spec(\Cox(X))$ d'\'equation
\begin{equation}
\prod_{i\in I'}s_i\neq 0
\end{equation}
pour $I'$ d\'ecrivant $\ecI_X$. Par ailleurs l'ouvert de $\Spec(\Cox(X))$ 
d'\'equation  $\produ{i\in I}s_i\neq 0$  est un ouvert non vide inclus
dans $\tors_X$. On le note $\tors_{X,0}$.
Pour tout $i\in I$, l'image r\'eciproque du support
du diviseur $\ecD_i$ est l'intersection du ferm\'e de $\Spec(\Cox(X))$
d'\'equation $s_i=0$ avec $\tors_X$.  En particulier
le ferm\'e d'\'equation $s_i=0$ rencontre $\tors_X$.
On a donc
\begin{equation}\label{eq:interIp}
\bigcap_{I'\in \ecI_X} I'=\vide.
\end{equation}
On a $X_0=\tors_{X,0}/\TNS$.

\begin{nota}\label{nota:tnx}
On note $T_{N_X}$ le tore $\Hom(N_X,\G_m)$.
La suite exacte \eqref{eq:suite_exacte_pic_tor} induit donc une suite exacte de tores
d\'eploy\'es
\begin{equation}\label{eq:suite_exacte_tor}
1\longto \TNS\longto 
\G_m^I \overset{\pi_X}\longto T_{N_X} \longto 1.
\end{equation}
\end{nota}
\begin{nota}\label{nota:nix}
On note $\nix $ le sous-ensemble de $\N^{\,I}$ form\'e des \'el\'ements $\bd$
v\'erifiant 
\begin{equation}
\forall \bn\in N_X,\quad \sum_{i\in I} n_i\,d_i=0.
\end{equation}
\end{nota}
\begin{rem}\label{rem:nix}
D'apr\`es la suite exacte \eqref{eq:suite_exacte_pic_tor} 
et la remarque \ref{rem:cone:effectif},
$\nix $ s'identifie \`a l'image dans $\Z^I$ de l'intersection de
$\ceff(X)$ avec $\Pic(X)^{\vee}$.
\end{rem}

\subsection{Application \`a la description du foncteur des points}
On se place toujours dans le cadre de la section
\ref{subsec:rappels:cox}, dont on conserve les notations.
Nous allons utiliser la description de $X$ comme le quotient g\'eom\'etrique $\tors_X/\TNS$ pour expliciter le foncteur des points
de $X$. Le r\'esultat est en fait une g\'en\'eralisation imm\'ediate et naturelle de la description de Cox du foncteur des points d'une
vari\'et\'e torique projective et lisse donn\'ee dans \cite{Cox:funct}.  
\begin{defi}\label{defi:nxtriv}
Soit $S$ un sch\'ema
et $(\ecL_i)_{i\in I}$ une famille de fibr\'es en droites 
sur $S$.
Une \termin{$N_X$-trivialisation} de $(\ecL_i)_{i\in I}$
est la donn\'ee d'une famille $(c_{\bn})_{\bn\in N_X}$ d'isomorphismes
\begin{equation}
c_{\bn}\,:\,\otimesu{i} \,\ecL_i^{n_i}\longisom \ecO_{S}
\end{equation}
v\'erifiant la conditions suivante :
pour tout $\bn,\bn'$ dans $N_X$, on a
\begin{equation}\label{eq:cond:cn}
c_{\bn}\otimes c_{\bn'}=c_{\bn+\bn'}.
\end{equation}
\end{defi}
\begin{rem}\label{rem:diff:cn}
Deux $N_X$-trivialisations $(c_{\bn})$ et $(c'_{\bn})$ de $(\ecL_i)_{i\in
  I}$ diff\`erent par 
une unique famille d'isomorphismes $c''_{\bn}\,:\,\ecO_S\isom\ecO_S$ v\'erifiant
la condition \eqref{eq:cond:cn},
{\it i.e.} par un morphisme de groupe
$
N_X\longto H^0(S,\ecO_S)^{\inv}.
$
En d'autres termes, l'ensemble des $N_X$-trivialisations de $(\ecL_i)_{i\in I}$
est un espace principal homog\`ene sous l'action du groupe $T_{N_X}\left(H^0(S,\ecO_S)\right)$.
\end{rem}
\begin{nota}\label{nota:fsi=0}
Pour tout \'el\'ement $D$ de $\Pic(X)$ on note
\begin{equation}
\N^{\,I}_{D}\eqdef\{\bd\in \N^{\,I},\quad \sumu{i} d_i\,D_i=D\}.
\end{equation}
Soit $F$ un polyn\^ome $\Pic(X)$-homog\`ene de $k[s_i]$,
de degr\'e $D$. \'Ecrivons
\begin{equation}
F=\sum_{\bd\in \N^{\,I}_{D}} \alpha_{\bd}\,\prod_{i\in I}s_i^{d_{i}}.
\end{equation}
Soit $S$ un $k$-sch\'ema.
Soit $(\ecL_i)_{i\in I}$ une famille de fibr\'es en droites sur $S$
munie d'une $N_X$-trivialisation $(c_{\bn})$. 
Pour tout $\ba=(a_i)\in \Z^I$ tel que 
$
D=\sumu{i\in I} a_{i}\,D_i
$
on pose
$
\ecL_{\ba}\eqdef\otimesu{i\in I}\ecL_i^{a_{i}}.
$ 
Pour tout $\bd\in \N^{\,I}_{D}$ tel que $\alpha_{\bd}\neq 0$
on a donc $\bd-\ba\in N_X$.
L'isomorphisme $c_{\bd-\ba}$  induit un isomorphisme
\begin{equation}
\varphi_{\bd,\ba}\,:\,\otimesu{i\in I}\ecL_i^{d_i}\longisom \ecL_{\ba}
\end{equation}
Si $\ba'=(a_i)$ est tel que
$
D=\sumu{i\in I} a'_{i}\,D_i
$, alors $\ba-\ba'\in N_X$
et l'isomorphisme $c_{\ba-\ba'}$ 
induit un isomorphisme
\begin{equation}
\iota_{\ba,\ba'}\,:\,\ecL_{\ba}\longisom \ecL_{\ba'}
\end{equation}
qui gr\^ace \`a la condition \eqref{eq:cond:cn} v\'erifie 
$\iota_{\ba,\ba'}\circ \varphi_{\bd,\ba}=\varphi_{\bd,\ba'}$
pour tout $\bd\in \N^{\,I}_{D}$.

Ainsi pour toute famille de sections globales $(u_i)$
on a 
\begin{equation}\label{eq:cond:annul}
\sum_{\bd\in \N^{\,I}_{D}} \alpha_{\bd}\,\varphi_{\bd,\ba}\left(\otimesu{i\in I} u_i^{\otimes d_{i}}\right)=0
\end{equation}
si et seulement si 
\begin{equation}
\sum_{\bd\in \N^{\,I}_{D}} \alpha_{\bd}\,\varphi_{\bd,\ba'}\left(\otimesu{i\in I} u_i^{\otimes d_{i}}\right)=0.
\end{equation}
La condition \eqref{eq:cond:annul} sera alors not\'ee
\begin{equation}
F(u_i)=0.
\end{equation}
On omet donc dans cette notation la r\'ef\'erence \`a la
$N_X$-trivialisation 
(qui dans la suite sera toujours clairement indiqu\'ee par le contexte).
\end{nota}

\begin{defi}\label{defi:collec}
Soit $S$ un $k$-sch\'ema.
Une \termin{$X$-collection sur $S$} est la donn\'ee pour tout 
$i\in I$ d'un fibr\'e en droites $\ecL_i$
sur $S$ et d'une section globale $u_i$ de $\ecL_i$,
ainsi que d'une $N_X$-trivialisation  $(c_{\bn})$ de $(\ecL_i)$,
ces donn\'ees \'etant astreintes \`a v\'erifier les conditions suivantes :
\begin{enumerate}
\item\label{item:defi:collec:prim}
Pour tout $i\in I$, la section 
$u_i$ induit un morphisme
$\ecO_S\to \ecL_i$
et par dualit\'e un morphisme
$\ecL_i^{-1}\to \ecO_S$.
On demande que le morphisme induit
\begin{equation}
\bigoplus_{I'\in \ecI_X} 
\,\underset{i\in I'}{\otimes}
\ecL_i^{-1}\longto \ecO_S
\end{equation}
soit surjectif.
\item\label{item:defi:collec:rel}
Pour tout \'el\'ement homog\`ene $F$ de $\idex$, on a
\begin{equation}
F(u_i)=0.
\end{equation}
\end{enumerate}
Un isomorphisme entre deux $X$-collections 
$\left((\ecL_i,u_i),(c_{\bn})\right)$ et
$\left((\ecL'_i,u'_i),(c'_{\bn})\right)$ 
est une famille d'isomorphismes $\ecL_i\isom \ecL'_i$
envoyant $u_i$ sur $u'_i$ et $c_{\bn}$ sur $c'_{\bn}$.
\end{defi}

On note $\cC_X(S)$ l'ensemble des $X$-collections sur $S$ modulo isomorphisme.
\begin{rem}
Il existe sur $X$ une $X$-collection universelle : 
on prend $\ecL_i=\ecO(\ecD_i)$ et $u_i$ la section canonique de
$\ecO(\ecD_i)$. 
La $N_X$-trivialisation $(c_{\bn})$
est donn\'ee par la suite exacte \eqref{eq:suite_exacte_pic_tor}
(elle d\'epend du choix de l'identification $N_X\subset k[X_0]^{\inv}$).
Si  $\pi\,:\,S\to X$ est un $k$-morphisme, 
$((\pi^{\ast}\ecO(\ecD_i),\pi^{\ast}u_i),(\pi^{\ast}c_{\bn}))$ est une
$X$-collection sur $S$. On obtient ainsi une application fonctorielle
en $S$
\begin{equation}\label{eq:homtocxs}
\Hom_{k}(S,X)\longto \cC_X(S).
\end{equation} 
\end{rem}
\begin{thm}\label{thm:cox}
L'application \eqref{eq:homtocxs} induit une bijection de 
$\Hom_{k}(S,X)$ sur l'ensemble des classes d'isomorphisme de $X$-collections sur $S$.
Ainsi $X$ repr\'esente
le foncteur qui \`a un $k$-sch\'ema $S$ associe l'ensemble des classes
d'isomorphisme de $X$-collections
sur $S$.
\end{thm}
La d\'emonstration est une adaptation imm\'ediate de la d\'emonstration du
th\'eor\`eme principal de \cite{Cox:funct}.
Indiquons juste de mani\`ere informelle comment est construit le morphisme
$S\to X$ correspondant \`a la $X$-collection
$\left((\ecL_i,u_i),(c_{\bn})\right)$ :
le morphisme en question associe \`a  $s\in S$ le 
\og point de coordonn\'ees
homog\`enes $(u_i(s))$\fg.
La suite exacte \eqref{eq:suite_exacte_pic_tor} montre
que le $I$-uple $(u_i(s))$ est bien d\'efini modulo l'action de $\TNS$.
Les conditions  
\ref{item:defi:collec:prim} 
et
\ref{item:defi:collec:rel} 
de la d\'efinition \ref{defi:collec}
assurent que  $(u_i(s))$ est dans $\ecT_X$.

On note $\Hom_{k,X_0}(S,X)$ l'ensemble des $k$-morphismes de $S$ vers
$X$ dont l'image sch\'ematique rencontre l'ouvert $X_0$.
\begin{defi}
Soit $S$ un $k$-sch\'ema.
Une $X$-collection sur $S$ 
$((\ecL_i,u_i),(c_n))$
est dite 
\termin{non d\'eg\'en\'er\'ee} 
si les sections $u_i$ sont toutes non nulles.
\end{defi}
Compte tenu du fait qu'on a $X_0=\tors_{X,0}/\TNS$ 
o\`u $\tors_{X,0}$ est l'ouvert d'\'equation $\prod s_i\neq 0$,
une adaptation imm\'ediate de la d\'emonstration du th\'eor\`eme \ref{thm:cox}
permet \'egalement de montrer le r\'esultat suivant. 
\begin{thm}\label{thm:nondegenerees}
L'application \eqref{eq:homtocxs} induit une bijection de 
$\Hom_{k,X_0}(S,X)$ sur l'ensemble des $X$-collections sur $S$ non
d\'eg\'en\'er\'ees modulo ismorphisme.
\end{thm}

\subsection{Description des morphismes de $\courbe$ vers $X$ : mont\'ee
  au torseur universel}
\label{subsec:courbe:vers:x}

On se place toujours dans le cadre de la section \ref{subsec:rappels:cox}.
Soit $\courbe$ une courbe projective, lisse et g\'eom\'etriquement int\`egre
d\'efinie sur $k$. On note $g$ son genre.
Pour $\bd\in \nix$, on cherche \`a d\'ecrire 
l'ensemble des $k$-morphismes de $\courbe$ vers $X$
dont l'image rencontre $X_0$ et tels qu'on ait 
\begin{equation}
\forall i\in I,\quad \deg f^{\ast} (D_i)=d_i.
\end{equation}
On va montrer que cet ensemble se d\'ecrit bien en termes
des donn\'ees d\'efinissant le torseur universel au dessus de $X$.

On suppose pour simplifier que $\courbe$ admet un diviseur de degr\'e
$1$
(ce qui sera de toute fa\c con v\'erifi\'e pour l'application que nous avons
en vue, le corps $k$ \'etant alors fini).
\begin{notas}
On {\em fixe} un diviseur de degr\'e $1$ sur $\courbe$, not\'e $\frD_1$.
On {\em fixe} \'egalement un sous-ensemble $\tpiczc $ de
$\Div^0(\courbe)$ de repr\'esentants de $\Pic^0(\courbe)$.

On note $\tpiczc^I_X$ le sous-ensemble de $\tpiczc ^I$ 
form\'e des \'el\'ements 
$(\frE_i)$ v\'erifiant 
\begin{equation}
\forall n\in N_X,\quad \sumu{i\in I} n_i\,\frE_i\sim 0.
\end{equation}
\end{notas}
\begin{rem}\label{rem:im:piczc}
L'image  de $\tpiczc^I_X$ dans $\Pic^0(\courbe)^I$
s'identifie donc au sous-groupe $\Pic^0(\courbe)\otimes N_X$.
\end{rem}
\begin{notas}
Pour tout $\bfrE\in \tpiczc^I_X$, 
on {\em fixe} une $N_X$-trivialisation $c_{\bfrE}$
de la famille $(\ecO_{\courbe}(\frE_i))$.
Pour tout $\bd\in \nix $, ceci induit
une $N_X$-trivialisation naturelle $c_{\bfrE,\bd}$ 
de la famille $(\ecO_{\courbe}(\frE_i+d_i\,\frD_1))$.
Par la suite, sauf mention explicite du contraire, 
la famille $(\ecO_{\courbe}(\frE_i+d_i\,\frD_1))$
sera toujours munie de cette trivialisation.

Pour tout $\frE\in \tpiczc $, 
et tout $d\in \N$
on note 
\begin{equation}
\Homog_{\frE,d}=H^0(\courbe,\ecO_{\courbe}(\frE+d\,\frD_1))
\end{equation}
et 
\begin{equation}
\Homogs_{\frE,d}=\Homog_{\frE,d}\setminus \{0\}.
\end{equation}
\end{notas}
Le lemme suivant 
est une cons\'equence classique du th\'eor\`eme de
Riemann-Roch et 
nous sera tr\`es utile lors de la d\'emonstration du r\'esultat
principal de cet article.
\begin{lemme}\label{lm:estim:dim}
\begin{enumerate}
\item
Si $d\geq 2\,g-1$, la dimension du $k$-espace vectoriel $\Homog_{\frD,d}$ est
$1-g+d$. 
\item
La dimension du $k$-espace vectoriel $\Homog_{\frD,d}$ est major\'ee par $1+d$.
\end{enumerate}
\end{lemme}

\begin{notas}\label{notas:nxd}
On note pour tout $\bfrE\in \tpiczc ^I$, 
et tout $\bd\in \N^{\,I}$
\begin{equation}
\Homogs_{\bfrE,\bd}
\eqdef
\prod_{i\in I}\Homogs_{\frE_i,d_i}.
\end{equation}
On note 
\begin{equation}
\Homogs_{I,X,\text{equiv}}
\eqdef
\disju{\substack{
\bfrE\in \tpiczc^I_X,\\ 
\bd\in \nix }}
\Homogs_{\bfrE,\bd}.
\end{equation}
\end{notas}

\begin{notas}
On note $\courbe^{(0)}$ l'ensemble des points ferm\'es de $\courbe$
et $\diveffc$ le mono\"\i de des diviseurs effectifs de $\courbe$,
{\it i.e.} le mono\"\i de ab\'elien libre de base $\courbe^{(0)}$.

Pour tout $v\in\courbe^{(0)}$ et $\ecD\in \diveffc$, on note $v(\ecD)$
la multiplicit\'e de $\ecD$ en $v$. 

Pour toute famille
$(\ecD_{\alpha})_{\alpha\in A}$ de diviseurs effectifs de $\courbe$,
on note
\begin{equation}
\pgcd\left((\ecD_{\alpha})_{\alpha\in A}\right)
\eqdef
\sum_{v\in \courbe^{(0)}} \Minu{\alpha\in A}\left(v(\ecD_{\alpha})\right)\,v.
\end{equation}
On pose
\begin{equation}
\Div_{I,X,\text{prim}}
\eqdef
\{\becE\in \diveffc^I,
\quad 
\underset{I'\in \ecI_X}{\pgcd}\left(\sum_{i\in I'}\ecE_i\right)=0
\}
\end{equation}
et
\begin{equation}
\Homogs_{I,X,\text{prim}}
\eqdef
\{(s_i)\in 
\!\!
\disju{
\substack{
\bfrE\in\tpiczc ^I,\\
\bd\in \N^{\,I}}
}
\Homogs_{\bfrE,\bd} ,
\quad
(\ddiv(s_i))\in \Div_{I,X,\text{prim}}
\}.
\end{equation}
\end{notas}

Pour tout $\bfrE\in \tpiczc^I_X$
et tout $\bd\in \nix $,
l'action diagonale naturelle de $\TNS(k)$ sur
$\Homogs_{\bfrE,\bd}$ 
pr\'eserve le sous-ensemble ({\it cf.} notations \ref{nota:fsi=0})
\begin{equation}
\{(s_i)\in \Homogs_{\bfrE,\bd} \cap \Homogs_{I,X,\text{prim}},\quad
\forall F\in \idexh, \quad F(s_i)=0\},
\end{equation}
o\`u $\idexh$ d\'esigne l'ensemble des \'el\'ements homog\`enes de l'id\'eal $\idex$.
\begin{prop}\label{prop:descr:mor}
Soit $\bd\in \nix $. 
On a une bijection  entre  :
\begin{enumerate}
\item
les $k$-morphismes de $\courbe$ vers $X$
dont l'image rencontre $X_0$ et v\'erifiant
\begin{equation}
\forall i\in I,\quad \deg f^{\ast} (\ecD_i)=d_i\quad ;
\end{equation}
\item\label{item:lm:descr:mor:2}
la r\'eunion des ensembles quotients
\begin{equation}
\{(u_i)\in \Homogs_{\frE,d} \cap \Homogs_{I,X,\text{prim}},\quad
\forall F\in \idexh, \quad F(s_i)=0\}/\TNS(k)
\end{equation}
pour $\bfrE$ parcourant $\tpiczc^I_X$.
\end{enumerate}
\end{prop}
\begin{demo}
\`A tout \'el\'ement de l'ensemble
\begin{equation}
\{(u_i)\in \Homogs_{\frE,d} \cap \Homogs_{I,X,\text{prim}},\quad
\forall F\in \idexh, \quad F(s_i)=0\}
\end{equation}
on associe la $X$-collection non deg\'en\'er\'ee
\begin{equation}
\left((\ecO_{\courbe}(\frE_i+d_i\,\frD_1),u_i),c_{\bfrE,\bd}\right)
\end{equation}
(le fait que $(u_i)$ soit dans $\Homogs_{I,X,\text{prim}}$ signifie exactement que la 
condition \ref{item:defi:collec:prim} de la d\'efinition \ref{defi:collec}
est satisfaite).

On va montrer que ceci induit une bijection de 
l'ensemble d\'ecrit dans le point 
\ref{item:lm:descr:mor:2}
de l'\'enonc\'e de la proposition
sur l'ensemble des 
classes d'isomorphisme de $X$-collections non 
d\'eg\'en\'er\'ees 
$((\ecL_i,u_i),(c_{\bn}))$
sur $\courbe$ v\'erifiant $\deg(\ecL_i)=d_i$,
ce qui donnera le r\'esultat d'apr\`es le th\'eor\`eme
\ref{thm:nondegenerees}.

Toute telle classe d'isomorphisme contient un \'el\'ement de la forme
\begin{equation}
((\ecO_{\courbe}(\frE_i+d_i\ecD_1),u_i),t.c_{\bfrE,\bd}), 
\end{equation}
o\`u 
$(\frE_i)\in \tpiczc^I_X$, $u_i\in \Homogs_{\frE_i,d_i}$,
et ({\it cf.} la remarque \ref{rem:diff:cn}) $t$ est un 
\'el\'ement de $T_{N_X}(H^0(\courbe,\ecO_{\courbe}))=T_{N_X}(k)$.

Deux $X$-collections 
\begin{equation}
((\ecO_{\courbe}(\frE_i+d_i\,\frD_1),u_i),t.c_{\bfrE,\bd})
\end{equation}
et 
\begin{equation}
((\ecO_{\courbe}(\frE_i+d_i\,\frD_1),u'_i),t'.c_{\bfrE,\bd})
\end{equation}
 sont isomorphes 
si et seulement s'il existe un \'el\'ement $(\lambda_i)\in
 (k^{\inv})^I=\G_m^I(k)$ tel que $u'_i=\lambda_i\,u_i$ et 
$t'=\pi_{X}(\lambda_i)\,t$ ({\it cf.} notation \ref{nota:tnx}). 

Ainsi dans toute classe d'isomorphisme, on peut
trouver un \'el\'ement de la forme $(\ecO_{\courbe}(\frE_i+d_i\,\frD_1),(u_i),c_{\bfrE,\bd})$.
Par ailleurs les $X$-collections
\begin{equation}
((\ecO_{\courbe}(\frE_i+d_i\,\frD_1),u_i),c_{\bfrE,\bd})
\end{equation}
et
\begin{equation}
((\ecO_{\courbe}(\frE_i+d_i\,\frD_1),u'_i),c_{\bfrE,\bd}) 
\end{equation}
sont isomorphes
si et seulement s'il existe un \'el\'ement $(\lambda_i)\in \G_m^I(k)$ tel que $u_i=\lambda_i\,u'_i$ et 
$\pi_X(\lambda_i)=1$, {\it i.e.} d'apr\`es \eqref{eq:suite_exacte_tor}
si et seulement s'il existe un \'el\'ement $(\lambda_i)\in \TNS(k)$ tel que $u_i=\lambda_i\,u'_i$. 
Ceci montre le r\'esultat.
\end{demo}

\subsection{Inversion de M\"obius}
On se place toujours dans le cadre de 
la section \ref{subsec:courbe:vers:x}.
Afin de se \og d\'ebarasser\fg\ de la condition 
$(s_i)\in \Homogs_{I,X,\text{prim}}$ apparaissant
dans la proposition \ref{prop:descr:mor}, il est classique
d'utiliser une inversion de M\"obius.

\begin{prop}\label{prop:mu}
Il existe une unique fonction $\mu_X\,:\,\diveffc^I\longto \C$
v\'erifiant
\begin{equation}
\forall \,\becD\in \diveffc^I,\quad
\ind_{\Div_{I,X,\text{prim}}}(\becD)
=
\sum_{0\leq \ecE_i\leq \ecD_i}
\mu_X(\becE)
\end{equation}
Cette fonction v\'erifie en outre les propri\'et\'es suivantes :
\begin{enumerate}
\item
 elle est multiplicative, c'est-\`a-dire que si $\becE$ et $\becD$ v\'erifient 
\begin{equation}
\forall i\in I,\quad \pgcd(\ecD_i,\ecE_i)=0
\end{equation}
alors on a
\begin{equation}
\mu_X(\becD+\becE)=\mu_X(\becD)\,\mu_X(\becE)\quad ;
\end{equation}
\item
pour tout $v\in\courbe^{(0)}$ et tout $\bn\in \N^{\,I}$,  
$\mu_{X}(\,(n_{i}\,v)\,)$ ne d\'epend que de $\bn$ (et pas de $v$) ; 
on note $\mu_{X}^0(\bn)$ cette valeur. 
On a $\mu_{X}^0(\bn)=0$ s'il existe $i$ tel que $n_{i}\geq 2$ ou  si
$\sum\,n_{i}=1$. 
\end{enumerate}
\end{prop}
\begin{demo}
Il suffit de reprendre la d\'emonstration de la proposition 1 de \cite{Bou:vtetor}.
Notons que la derni\`ere assertion d\'ecoule aussit\^ot du fait
que si $\sum\,n_{i}=1$ alors 
$
((n_i\,v))$ est dans  $\Div_{I,X,\text{prim}}$
ce qui provient de $\eqref{eq:interIp}$.
\end{demo}
\begin{rem}\label{rq:mu0}
Notons $\{0,1\}^I_X$ l'ensemble des \'el\'ements $(n_i)\in \{0,1\}^I$
v\'erifiant
\begin{equation}
\Minu{I'\in \ecI_X} \,\,\sum_{i\in I'} n_i=0
\end{equation} 
On a alors
\begin{equation}\label{eq:ind01IX}
\forall \bn\in \{0,1\}^I,\quad \ind_{\{0,1\}^I_X}(\bn)=\sum_{0\leq \bm
  \leq \bn} \mu^0_X(\bm).
\end{equation}
En particulier, si $\bn\in \{0,1\}^I_X\setminus\{0\}$, on a
$\mu^0_X(\bn)=0$.
\end{rem}
\begin{nota}\label{nota:LA}
Si $A$ est un ensemble, $\bn$ un \'el\'ement de $\{0,1\}^A$
et $L$ un corps, on pose
\begin{equation}
L^{\bn}\eqdef
\{(x_{\alpha})\in L^{A},\quad
\forall \alpha\in A,\quad x_{\alpha}=0\,\text{ si }\,n_{\alpha}=1
\}.
\end{equation}
\end{nota}

On suppose \`a pr\'esent que $k$ est un corps fini de cardinal $q$.
Pour $v\in \courbe^{(0)}$, on note $f_v$ le degr\'e de $v$ et $\kappa_v$
le corps r\'esiduel, de sorte que $\card{\kappa_v}=q^{f_v}\eqdef q_v$.
\begin{nota}\label{nota:dens}
Pour $\bn\in \{0,1\}^I$,
on pose
\begin{equation}
\dens_{X,v}(\bn)
\eqdef
\frac{
\card{
\{(x_i)\in \kappa_v^{\bn},\quad \forall F\in \idex,\quad F(x_i)=0\}
}}
{
q_v^{\dim(\ecT_X)}.
}
\end{equation}
\end{nota}
\begin{lemme}\label{lm:rel:densv}
On a la relation
\begin{equation}
\sum_{\bn\in \{0,1\}^I}\,\mu^0_X(\bn)\,\dens_{X,v}(\bn)
=
(1-q_v^{-1})^{\rg(\TNS)}
\frac{\card{X(\kappa_v)}}
{q_v^{\dim(X)}}.
\end{equation}
\end{lemme}
\begin{demo}
On a 
\begin{equation}
\forall (x_i)\in \kappa_v^I, 
\left(\exists I'\in \ecI_X, \prod_{i\in I'}x_i\neq 0
\longeq
\forall \bn\notin \{0,1\}^I_X,\quad (x_i)\notin \kappa_v^{\bn}
\right)
\end{equation}
On en d\'eduit l'\'egalit\'e 
\begin{equation}
\tors_X(k_v)
=
\left\{(x_i)\in \kappa_v^I, \quad \forall \bn\notin
\{0,1\}_I^X,\,(x_i)\notin \kappa_v^{\bn},\quad
\forall F\in \idex,\,F(x_i)=0
\right\}
\end{equation}
Ainsi, d'apr\`es \eqref{eq:ind01IX}, on a
\begin{equation}
\ind_{\tors_X(k_v)}
=
\sum_{\bn\in \{0,1\}^I}\,\mu^0(\bn)\,\ind_{\{(x_i)\in k_v^{\bn}, \quad
\forall F\in \idex,\, F(x_i)=0\}}.
\end{equation}
On en tire
\begin{equation}
\frac
{\card{\tors_X(\kappa_v)}}
{q_v^{\dim(\tors_X)}}
=
\sum_{\bn\in \{0,1\}^I}\,\mu^0(\bn)\,\dens_{X,v}(\bn)
\end{equation}
Pour conclure il suffit de remarquer que comme $X$
est le quotient g\'eom\'etrique $\ecT_X/\TNS$
et  $\TNS$ est d\'eploy\'e on a
\begin{equation}
\card{X(\kappa_v)}=\frac{\card{\tors_X(\kappa_v)}}{(q_v-1)^{\rg(\TNS)}}
\end{equation}
et $\dim(\tors_X)=\rg(\TNS)+\dim(X)$.
\end{demo}

\subsection{Conjectures de Manin}
On se place toujours dans le cadre de 
la section \ref{subsec:courbe:vers:x}.
On suppose d\'esormais pour tout le reste de l'article que le corps de
base $k$ est fini
de cardinal $q$.
On note
\begin{equation}
Z_{\courbe}(T)\eqdef\sum_{\ecD\in \diveffc}T^{\deg(\ecD)}
\end{equation}
la fonction z\^eta de Dedekind de la courbe $\courbe$.

Soit $D$ un \'el\'ement de $\Pic(X)$ situ\'e \`a l'int\'erieur du c\^one effectif.
Pour tout ouvert $U$ de $X$, on note, pour $d\geq 0$,
$
N_{D,U}(d)
$
le cardinal  de l'ensemble des $k$-morphismes $f\,:\,\courbe\to X$
dont l'image recontre $U$ et qui v\'erifient
\begin{equation}
\deg f^{\ast}(D)=d.
\end{equation}
Ce cardinal est fini si $U$ est assez petit (et toujours fini si $D$
est ample par exemple).

Les conjectures de Manin tentent alors de pr\'edire le comportement
asymptotique de $N_{D,U}(d)$ lorsque $d$ tend vers l'infini.
Nous \'enon\c cons une version possible de cette conjecture dans le cas o\`u
$D$ est la classe du fibr\'e anticanonique.

Pour cela, on d\'efinit, suivant Peyre, les constantes 
\begin{equation}
\alpha(X)
\eqdef
\lim_{T\to 1} (1-T)^{\rg(\Pic(X))}\,
\sum_{y\in \ceff(X)^{\vee}\cap \Pic(X)^{\vee}}\,T^{\,\acc{y}{\classe{\omega_X^{-1}}}}
\end{equation}
et
\begin{equation}
\gamma(X)
\eqdef
\left(\lim_{T\to q^{-1}}
    (1-q\,T)\,Z_{\courbe}(T)\right)^{\rg(\Pic(X))}\:
q^{(1-g)\,\dim(X)}\,
\prod_{v\in\courbe^{(0)}}
(1-q_v^{-1})^{\rg(\Pic(X))}\,\frac{\card{X(\kappa_{v})}}{q_v^{\,\dim(X)}}
\end{equation}
Nous renvoyons \`a \cite{Pey:var_drap} pour la justification du fait que ces constantes
sont bien d\'efinies. L'argument de loin le plus d\'elicat concerne la
convergence du produit eul\'erien figurant dans la d\'efinition de
$\gamma(X)$, pour lequel on invoque les conjectures de Weil d\'emontr\'ees par Deligne.
Il est \`a noter que dans tous les cas o\`u la conjecture de Manin a \'et\'e \'etablie, la
convergence peut \^etre d\'emontr\'ee directement, sans faire appel \`a un
r\'esultat aussi fin.

\begin{question}\label{ques:manin:1}
Soit 
\begin{equation}
\delta=\Max \{d\in \N_{>0},\quad \frac{1}{d}\classe{\omega_X^{-1}}\in \Pic(X)\}.
\end{equation}
A-t-on, pour tout ouvert $U$ de $X$ assez petit,
\begin{equation}
N_{\omega_X^{-1},U}(\delta\,d)\simu{d\to +\infty}
\alpha(X)\,\gamma(X)\,d^{\,\,\rg(\Pic(X))-1} q^{\,\delta\,d}\quad ?
\end{equation}
\end{question}
Une strat\'egie classique pour l'\'etude asymptotique de 
$
N_{D,U}(d)
$
est d'essayer de pr\'eciser le comportement analytique 
de la fonction z\^eta des hauteurs associ\'ee, {\it i.e.} la s\'erie g\'en\'eratrice 
\begin{equation}
Z_{D,U}(T)\eqdef \sum_{d\geq 0} N_{D,U}(d)\,T^d.
\end{equation}
Dans cette optique, rappelons d'abord deux \'enonc\'es taub\'eriens
\'el\'ementaires, cons\'equences directes des estimations de Cauchy.
\begin{prop}\label{prop:tauber}
Soit $(a_n)\in \C^{\N}$,
$\alpha\in \C^{\ast}$
et $k\geq 1$ un entier. 
On suppose que la s\'erie $\sum a_n\,z^n$ a pour rayon de convergence
$\abs{\alpha}$ et que sa somme se prolonge en une fonction $f(z)$ m\'eromorphe
sur un disque de rayon strictement sup\'erieur \`a $\abs{\alpha}$, ayant en
$\alpha$
un p\^ole d'ordre $k$ et des p\^oles d'ordre au plus $k-1$ en tout autre
point
du cercle de rayon $\abs{\alpha}$.
Alors on a
\begin{equation}
a_n=(\lim_{z\to \alpha}(z-\alpha)^k\,f(z))\,n^{k-1}\,\alpha^{-n}
+
\ecOu{n\to +\infty}\left(n^{k-2}\,\abs{\alpha}^{-n}\right).
\end{equation}
\end{prop}
\begin{defi}
On dit que la s\'erie $\sum a_n\,z^n$ 
(\`a coefficients complexes)
 est major\'ee par par la s\'erie $\sum b_n\,z^n$ 
(\`a coefficients r\'eels positifs) si on a $\abs{a_n}\leq b_n$ pour
tout $n$.
\end{defi}
\begin{prop}\label{prop:tauber:2}
Soit $(a_n)\in \C^{\N}$, $k\geq 1$ un entier et $\rho>0$ un r\'eel.
Les conditions suivantes sont \'equivalentes :
\begin{enumerate}
\item
on a 
\begin{equation}
a_n=\ecOu{n\to +\infty}\left(n^{k-1}\,\rho^{-n}\right) ;
\end{equation}
\item
la s\'erie $\sum a_nz^n$ est major\'ee par une s\'erie 
dont le rayon de convergence est sup\'erieur \`a $\rho$
et dont la somme se prolonge en une fonction m\'eromorphe
sur un disque de rayon strictement sup\'erieur  \`a $\rho$,
ayant des p\^oles d'ordre au plus $k$ sur le cercle de rayon $\rho$.
\end{enumerate}
\end{prop}
\begin{defi}
On dit que $\sum a_nz^n$ est $\rho$-contr\^ol\'ee \`a l'ordre $k$
si elle v\'erifie les conditions de la proposition \ref{prop:tauber:2}.
\end{defi}

Dans le cas o\`u $D=\classe{\omega_X^{-1}}$ on peut \'enoncer
la variante analytique suivante de la question \ref{ques:manin:1}.
\begin{question}\label{ques:manin:2}
On conserve les notations de la question \ref{ques:manin:1}.
On note $\wt{Z}_{\omega_{X}^{-1},U}(T)$ la s\'erie telle que 
$\wt{Z}_{\omega_{X}^{-1},U}(T^{\delta})=Z_{\omega_{X}^{-1},U}(T)$.
Est-il vrai que si $U$ est assez petit, 
la s\'erie $\wt{Z}_{\omega_{X}^{-1},U}(T)$
a pour rayon de convergence  $q^{-\delta}$ et que 
sa somme 
se prolonge en une
fonction m\'eromorphe sur le disque $\abs{z}<q^{-\delta+\eps}$ ayant un 
p\^ole d'ordre $\rg(\Pic(X))$ en $z=q^{-\delta}$, 
et des p\^oles d'ordre
au plus $\rg(\Pic(X))-1$ en tout autre point du cercle de
rayon $q^{-\delta}$, et v\'erifiant
\begin{equation}
\lim_{T\to q^{-\delta}}\left(T-q^{\,-\delta}\right)^{\rg(\Pic(X))}\,
\wt{Z}_{\omega_{X}^{-1},U}(T)=\alpha(X)\,\gamma(X)
\end{equation}
\end{question}
Si la question \ref{ques:manin:1} admet une r\'eponse positive, 
alors d'apr\`es la proposition \ref{prop:tauber} la question \ref{ques:manin:2} admet une r\'eponse positive.

\subsection{Rel\`evement au torseur universel pour la fonction z\^eta des hauteurs}\label{subsec:montee:tu:zeta}

On consid\`ere toujours un \'el\'ement $D$ de $\Pic(X)$
situ\'e \`a l'int\'erieur du c\^one effectif.
Nous nous pla\c cons \`a pr\'esent dans le cas o\`u $U=X_0$
et expliquons comment s'exprime au niveau de la fonction z\^eta des
hauteurs la mont\'ee au torseur universel (donn\'ee par la proposition
\ref{prop:descr:mor})
du d\'ecompte des morphismes
de degr\'e born\'e.

Soit $(n_{i,D})\in \N_{>0}^I$ tel que 
$
D=\sum_i n_{i,D}\,\classe{\ecD_i}.
$
Nous avons d'apr\`es les propositions
\ref{prop:descr:mor}
et \ref{prop:mu}
\begin{align}
(q-1)^{\rg(\TNS)}\,Z_{D,X_0}(T)
&=
\sum_{
\substack{
(s_i)\in \Homogs_{I,X,\text{equiv}}\cap \Homogs_{I,X,\text{prim}}
\\ \\
\forall F\in \idexh,\quad F(s_i)=0
}
}
\,T^{\,\sumu{i}n_{i,D} \deg(s_i)} 
\\  
&=\sum_{
\substack{
(s_i)\in \Homogs_{I,X,\text{equiv}}
\\ \\
\forall F\in \idexh,\quad F(s_i)=0
}
}\,
\left(
\sum_{
\substack{
\becE\in \diveffc^{I} \\ \\
 \ecE_i\leq \ddiv(s_i)}} 
\mu_X(\becE)
\right) 
\, 
T^{\,\sumu{i}n_{i,D} \deg(s_i)} 
\\
&=
\sum_{\becE\in \diveffc^{I}}
\,
\mu_X(\becE)
\,
\left(
\sum_{
\substack{ 
(s_i)\in {\Homogs_{I,X,\text{equiv}}}
\\ \\
\ddiv(s_i)\geq \ecE_i
\\ \\
\forall F\in \idexh,\quad F(s_i)=0
}
} 
\, 
T^{\,\sumu{i}n_{i,D} \deg(s_i)}.
\right) \label{eq:egalite}
\end{align}
\begin{nota}
Pour tout \'el\'ement $\ecE\in \diveffc$, on note $\scan{\ecE}$ 
la section canonique de $\ecO_{\courbe}(\ecE)$.
\end{nota}

\begin{notas}\label{nota:piczcixe}
Soit $\becE\in \diveffc^I$. 
On note $\piczcixe$ 
l'ensemble des \'el\'ements $\bfrE\in\tpiczc ^I$ v\'erifiant
\begin{equation}
\bfrE+(\ecE_i-\deg(\ecE_i)\,\frD_1) \in \tpiczc^I_X.
\end{equation} 
Remarquons que le morphisme \og classe dans le groupe de Picard\fg\
induit une bijection de $\piczcixe$ sur le sous-ensemble de $\Pic^0(\courbe)^I$
donn\'e par
\begin{equation}
-(\classe{\ecE_i}-\deg(\ecE_i)\,\classe{\frD_1})+\Pic^0(\courbe)\otimes N_X.
\end{equation}

Pour $\bd\in \nix$ v\'erifiant $d_i\geq \deg(\ecE_i)$,
et $\bfrE\in\tpiczc ^I$,
on note $\ecN_X(\bd,\becE,\bfrE)$ le cardinal de l'ensemble des \'el\'ements
\begin{equation}
(s_i)\in \prod_{i\in I}\Homogs_{\frE_i,d_i-\deg(\ecE_i)}
\end{equation}
v\'erifiant
\begin{equation}
\forall F\in \idexh,\quad F(s_i.\scan{\ecE_i})=0.
\end{equation}
\end{notas}

Ainsi on a
\begin{align}
&\phantom{=}(q-1)^{\rg(\TNS)}\,
Z_{D,X_0}(T)
\notag
\\
&
=
\sum_{\becE\in \diveffc^{I}}
\,
\mu_X(\becE)
\sum_{\bfrE\in
\piczcixe }
\left(
\sum_{\substack{\bd\in \nix , \\ d_i\geq \deg(\ecE_i)}}
\ecN_X(\bd,\becE,\bfrE)
\, 
T^{\,\sumu{i}n_{i,D} d_i}
\right)
\label{eq:expr:zdx0}
\end{align}

\begin{rem}
La remarque \ref{rem:nix} permet de r\'e\'ecrire 
l'\'egalit\'e \eqref{eq:expr:zdx0}
en termes du c\^one effectif de $X$, {\it i.e.}
sous la forme
\begin{align}
&\phantom{=}(q-1)^{\rg(\TNS)}\,
Z_{D,X_0}(T)
\notag
\\
&
=
\sum_{\becE\in \diveffc^{I}}
\,
\mu_X(\becE)
\sum_{\bfrE\in
\piczcixe }
\left(
\sum_{
\substack{ 
y\in \Pic(X)^{\vee}\cap \ceff(X)^{\vee}
\\ \\
\acc{y}{D_i}\geq \deg(\ecE_i)
}
}
\ecN_X((\acc{y}{D_i}),\becE,\bfrE)
\, 
T^{\,\,\acc{y}{D}}
\right) .
\end{align}
\end{rem}
Comme dans le cadre des corps de nombres, 
cette mont\'ee au torseur universel ne constitue qu'une
premi\`ere \'etape dans une \'eventuelle d\'emonstration de la conjecture de
Manin pour $X$. La t\^ache difficile consiste \`a \'evaluer de mani\`ere
suffisamment
pr\'ecise le comportement asymptotique de la quantit\'e $\ecN_X(\bd,\becE,\bfrE)$.
Comme d\'ej\`a indiqu\'e dans l'introduction,  le cas le plus favorable 
\`a cet \'egard est celui o\`u $X$ est torique, car l'id\'eal $\idex$ est alors
nul : il n'y a pas d'\'equation \`a prendre en compte.
Dans la suite de cet article, nous expliquons comment traiter le cas
du plan projectif \'eclat\'e en trois points align\'es, cas o\`u le torseur
universel est donn\'e par une unique \'equation, qui plus est
particuli\`erement simple.

\section{Le cas du plan projectif \'eclat\'e en 3 points align\'es}

\subsection{Description du torseur universel au-dessus du plan \'eclat\'e en 3 points align\'es}
\label{subsec:tors:p2}
Soit $p_1$, $p_2$, $p_3$ trois points align\'es du plan projectif et $S$
la surface obtenue en \'eclatant ces trois points.
Soit $p_0$ un point du plan projectif qui n'est pas sur la droite
$\langle p_1,p_2,p_3\rangle$.
On note $(\cE_i)_{i=1,2,3}$ les diviseurs exceptionnels de
l'\'eclatement, $\cE_0$ le transform\'e strict de la droite
$\langle p_1,p_2,p_3 \rangle$ et $(\cF_i)_{i=1,2,3}$ les transform\'es stricts
de droites $\langle p_i,p_0\rangle$.

D'apr\`es \cite{Der:sdp:ut:hyp} ({\it cf.} \'egalement
\cite{Has:eq:ut:cox:rings}),
on peut trouver des sections globales $s_i$ (respectivement $t_i$)
de diviseur $\cE_i$ (respectivement $\cF_i$) tel qu'on ait un isomorphisme
\begin{equation}
\Cox(S)\longisom k[(s_i)_{i=0,\dots,3},(t_i)_{i=1,2,3}]/
\sum_{i=1}^3 s_i\,t_i
\end{equation}
On a 
\begin{equation}
\Pic(S)=\oplusu{0\leq i \leq 3} \Z\,\classe{\cE_i}
\end{equation}
et pour $i=1,2,3$
\begin{equation}
\classe{\cF_i}=\sum_{j=0,j\neq i}^3\classe{\cE_i}.
\end{equation}

La classe du fibr\'e anticanonique
est
\begin{equation}
\classe{\omega_S^{-1}}=3\,\classe{\cE_0}+2\sum_{i=
1}^3\classe{\cE_i}.
\end{equation}

On v\'erifie alors que
$\tors_S\subset \Spec(\Cox(S))$ est la r\'eunion des
ouverts d'\'equations
\begin{align}
s_1\,s_2\,t_1\,t_2\,t_3&\neq 0,
\\
s_2\,s_3\,t_1\,t_2\,t_3&\neq 0,
\\
s_1\,s_3\,t_1\,t_2\,t_3&\neq 0,
\\
s_0\,s_1\,s_2\,t_1\,t_2&\neq 0,
\\
s_0\,s_1\,s_3\,t_1\,t_3&\neq 0,
\\
s_0\,s_2\,s_3\,t_2\,t_3&\neq 0,
\\
\text{et}\quad s_0\,s_1\,s_2\,s_3&\neq 0.
\end{align}
\begin{rem}\label{rq:ceffS}
On a d'apr\`es ce qui pr\'ec\`ede
\begin{equation}
\sum_{y\in \Pic(S)^{\vee}\cap \ceff(S)^{\vee}}
T^{\acc{y}{\classe{\omega_S^{-1}}}}=\frac{1}{(1-T^3)\,(1-T^2)^3}.
\end{equation}
\end{rem}
\begin{rem}\label{rq:nis}
D'apr\`es ce qui pr\'ec\`ede, l'application
\begin{equation}
(d_0,d_1,d_2,d_3)\longto (d_0,d_1,d_2,d_3,d_0+d_1+d_3,d_0+d_2+d_3,d_0+d_1+d_2)
\end{equation}
est une bijection de $\N^4$ sur $\N^7_S$ ({\it cf.} la notation \ref{nota:nix}). Par la suite, on identifiera 
toujours $\N^7_S$ \`a $\N^4$ au moyen de cette bijection.

De m\^eme, pour tout $\becE=(\ecE_0,(\ecE_i),(\ecF_i))
\in\diveffc^{7}$, l'application qui 
\`a $(D_0,D_1,D_2,D_3)\in \Pic^0(\courbe)^4$
associe
\begin{equation}
\left(D_0,(D_i),
\left(
D_0+\classe{\ecE_0}-\deg(\ecE_0)\classe{\frD_1}
+\sum_{j\neq i}\left(D_i+\classe{\ecE_i}-\deg(\ecE_i)\classe{\frD_1}\right)
-\classe{\ecF_i}+\deg(\ecF_i)\,\classe{\frD_1}
\right)
\right)
\end{equation}
est une bijection de $\Pic^0(\courbe)^4$ sur le sous-ensemble
de $\Pic^0(\courbe)^7$ donn\'e par
\begin{equation}
-\Big(
\ecE_0-\deg(\ecE_0)\,\classe{\frD_1},
(\ecE_i-\deg(\ecE_i)\,\classe{\frD_1}),
(\ecF_i-\deg(\ecF_i)\,\classe{\frD_1})
\Big)
+
\Pic^0(\courbe)\otimes N_S.
\end{equation}
Ainsi ({\it cf.} la notation \ref{nota:piczcixe}
et la remarque y figurant) l'application 
\begin{equation}
(\frE_0,(\frE_i),(\frF_i))\mapsto (\frE_0,\frE_1,\frE_2,\frE_3)
\end{equation}
compos\'ee avec le morphisme \og classe dans le groupe de Picard\fg\
induit une bijection de $\piczcSe$ sur $\Pic^0(\courbe)^4$.
\end{rem}

\subsection{Le r\'esultat}
\begin{thm}\label{thm:resprinc}
Soit $p_1$, $p_2$, $p_3$ trois points align\'es du plan projectif et $S$
la surface obtenue en \'eclatant ces trois points. Soit $p_0$
un point du plan projectif qui n'est pas sur la droite
$\langle p_1,p_2,p_3\rangle $.
Soit $S_0$ l'ouvert de $S$ obtenu en retirant les diviseurs exceptionnels de l'\'eclatement,
et les transform\'es stricts des droites
$\langle p_1,p_2,p_3 \rangle$, et $\langle p_i,p_0\rangle$ pour $i=1,2,3$.

On a alors une \'ecriture
\begin{equation}
Z_{\omega_{S}^{-1},S_0}(T)=
Z_{\courbe}(q^2\,T^3)\,Z_{\courbe}(q\,T^2)^3
\wt{Z}(T)
+Z_{\text{err}}(T)
\end{equation}
o\`u $\wt{Z}(T)$ est une s\'erie de rayon de convergence strictement
sup\'erieure \`a $q^{-1}$ v\'erifiant
\begin{equation}
\wt{Z}(q^{-1})=q^{\,2(1-g)}
\prod_{v\in\courbe^{(0)}}
(1-q_v^{-1})^{\rg(\Pic(S))}\,\frac{\card{S(\kappa_{v})}}{q_v^{\,\dim(S)}}
\end{equation}
et
$Z_{\text{err}}(T)$ est une s\'erie $q^{-1}$-contr\^ol\'ee \`a l'ordre $2$.
\end{thm}
La d\'emonstration de ce th\'eor\`eme fait l'objet du reste de cet article.
\begin{cor}
On a 
\begin{equation}
N_{\omega_S^{-1},S_0}(d)=
\alpha(S)\,\gamma(S)\,d^{\,\,\rg(\Pic(S))-1} q^{\,d}+\ecOu{d\to +\infty}\left(d^{\,\,\rg(\Pic(S))-2}\,q^{\,d}\right).
\end{equation}
En particulier, la r\'eponse \`a la question \ref{ques:manin:1} est positive pour l'ouvert $U=S_0$.
\end{cor}
\begin{demo}
Notons qu'on a ici $\delta=1$.
Compte tenu du th\'eor\`eme \ref{thm:resprinc}
et des propositions \ref{prop:tauber} et
 \ref{prop:tauber:2}, il suffit de montrer qu'on a 
\begin{equation}
\lim_{T\to q^{-1}} \left(T-q^{-1}\right)^4 Z_{\courbe}(q^2\,T^3)\,Z_{\courbe}(q\,T^2)^3
=\alpha(S)\,\left(\lim_{T\to q^{-1}} (1-q\,T)\,Z_{\courbe}(T)\right)^4
\end{equation}
ce qui est imm\'ediat au vu de la remarque \ref{rq:ceffS}.
\end{demo}

\subsection{D\'emonstration du r\'esultat principal : pr\'eliminaires}\label{subsec:prelim}
Rappelons ({\it cf.} la remarque \ref{rq:nis}) 
que l'on a identifi\'e $\N^7_{S}$ \`a $\N^4$.
Ainsi, en reprenant la notation \ref{notas:nxd},
pour $\becE=(\ecE_0,(\ecE_i),(\ecF_i)) \in\diveffc^{7}$,
$\bfrE=(\frE_0,(\frE_i),(\frF_i))\in \piczcSe$,
et $\bd\in \N^4$,
$\ecN_S(\bd ,\becE ,\bfrE)$ d\'esigne
le cardinal de l'ensemble des \'el\'ements 
\begin{equation}
(s_0,(s_i),(t_i))\in \Homogs_{\frE_0,d_0-\deg(\ecE_0)}
\times \prod_i \Homogs_{\frE_i,d_i-\deg(\ecE_i)}
\times \prod_i \Homogs_{\frF_i,d_0+\sumu{j\neq i}d_j-\deg(\ecF_i)}
\end{equation}
v\'erifiant la relation 
\begin{equation}
\sumu{1\leq i\leq 3} s_i\,t_i\,\scan{\ecE_i}\,\scan{\ecF_i}=0.
\end{equation}
On pose
\begin{equation}\label{eq:exprzece}
Z(\becE ,\bfrE,T)
\eqdef
\sum_{
\substack{
d_0\geq \deg(\ecE_0) \\
d_i\geq \deg(\ecE_i) \\
d_0+\sumu{j\neq i}d_j\geq \deg(\ecF_i)
}
}
\ecN_S(\bd ,\becE ,\bfrE)
\,
T^{3\,d_0+2\,\sumu{i}d_i}.
\end{equation}

D'apr\`es les sections \ref{subsec:montee:tu:zeta} 
et \ref{subsec:tors:p2},
la fonction z\^eta des hauteurs $Z_{\omega_{S}^{-1},S_0}(T)$
s'\'ecrit alors
\begin{equation}
\frac{1}{(q-1)^4}
\sum_{\becE \in\diveffc^{7}}
\mu_S(\becE) 
\sum_{\bfrE\in \piczcSe }
Z(\becE ,\bfrE,T).
\end{equation}

On pose pour $1\leq i\leq 3$
\begin{equation}\label{eqref:psi}
\psi_i(\bd,\becE)
\eqdef 
d_0+\sumu{j\neq i}d_j +\deg(\ecE_0)+\sumu{j\neq i}\deg(\ecE_j)-\deg(\ecF_i).
\end{equation}
Afin d'all\'eger un peu l'\'ecriture, on change l\'eg\`erement les notations :
on d\'esigne d\'esormais par $\ecN_S(\bd ,\becE ,\bfrE)$
le cardinal de l'ensemble des \'el\'ements 
\begin{equation}
(s_0,(s_i),(t_i))\in 
\Homogs_{\frE_0,d_0} 
\times
\prod_i \Homogs_{\frE_i,d_i} 
\times
\prod_i \Homogs_{\frF_i,d_0+\sumu{j\neq i}d_j+\deg(\ecE_0)+\sumu{j\neq i}\deg(\ecE_j)-\deg(\ecF_i)}
\end{equation}
v\'erifiant la relation $\sumu{i} s_i\,t_i\,\scan{\ecE_i}\,\scan{\ecF_i}=0$.

Un changement de variables imm\'ediat dans l'expression \eqref{eq:exprzece}
permet donc d'\'ecrire  $Z(\becE ,\bfrE,T)$
sous la forme
\begin{equation}\label{eq:def}
T^{3\,\deg(\ecE_0)+2\,\sumu{i}\deg(\ecE_i)}
\sum_{
\substack{
\bd\in \N^4\\
~
\\
\forall 1\leq i\leq 3,\quad \psi_i(\bd ,\becE )\geq 0
}}
\ecN_S\left(
\bd,
\becE ,\bfrE
\right)
\,T^{3\,d_0+2\,\sumu{i}d_i}
\end{equation}

On d\'efinit
\begin{equation}
\phi_1(\bd ,\becE )\eqdef d_0+d_1+\deg(\ecE_0)+\deg(\ecE_1)-\deg(\ecF_2)-\deg(\ecF_3)
\end{equation}
et $\phi_2$, $\phi_3$ de mani\`ere analogue en permutant de mani\`ere
circulaire les indices $1$, $2$ et $3$.

Soit $\becE\in \diveffc^7$ et $\bfrE\in \piczcSe$.
Soit $\Zp(\becE,\bfrE,T)$ la s\'erie d\'efinie par la formule
\eqref{eq:def}
en restreignant le domaine de sommation aux $\bd $ v\'erifiant la
contrainte suivante : il existe deux indices distincts $k$ et $k'$ avec $1\leq k,k'\leq
3$ tels qu'on ait
\begin{equation}
\phi_{k}(\bd ,\becE )\geq 2\,g-1
\end{equation}
et
\begin{equation}
\phi_{k'}(\bd ,\becE )\geq 2\,g-1.
\end{equation}

Pour $1\leq k \leq 3$, on \'ecrit $\{1,2,3\}=\{k,k',k''\}$.
Soit $Z_k(\becE,\bfrE,T)$ la s\'erie d\'efinie par la formule
\eqref{eq:def}
en restreignant le domaine de sommation aux $\bd $ v\'erifiant les contraintes  
\begin{equation}
\phi_{k}(\bd ,\becE )\geq 2\,g-1,
\end{equation}
\begin{equation}
\phi_{k'}(\bd ,\becE )<2\,g-1
\end{equation}
et
\begin{equation}
\phi_{k''}(\bd ,\becE )<2\,g-1.
\end{equation}

Soit enfin $Z_4(\becE,\bfrE,T)$ la s\'erie d\'efinie par la formule
\eqref{eq:def}
en restreignant le domaine de sommation aux $\bd $ v\'erifiant les contraintes  
\begin{equation}
\forall i\in \{1,2,3\},\quad \phi_{i}(\bd ,\becE )<2\,g-1.
\end{equation}

On a donc l'\'ecriture
\begin{equation}
Z(\becE,\bfrE,T)=\Zp(\becE ,\bfrE,T)
+\sum_{k=1}^3 Z_{k}(\becE ,\bfrE,T)
+Z_{4}(\becE ,\bfrE,T).
\end{equation}

Expliquons en deux mots l'int\'er\^et de cette d\'ecomposition. Comme on le
verra \`a la section \ref{subsec:compt:sec}, nous allons utiliser le
th\'eor\`eme de Riemann-Roch pour estimer la quantit\'e $\ecN_S(\bd ,\becE
,\bfrE)$. Pour $\bd$ \og grand \fg, on obtiendra une formule exacte
alors que  pour $\bd$ \og petit \fg, on devra se contenter d'une majoration ({\it
  cf.} le corollaire \ref{cor:nsO}). Le terme $\Zp$ (respectivement les termes
$Z_1,\,\dots,\,Z_4$) correspond \`a la
sommation sur les $\bd$ qui sont \og grands\fg\ (respectivement \og
petits \fg).
Les majorations  de $\ecN_S(\bd ,\becE,\bfrE)$ obtenues \`a la section \ref{subsec:compt:sec}
permettront de montrer que les termes $Z_1,\,\dots,\,Z_4$ ne contribuent
pas au terme principal de la fonction z\^eta des hauteurs. Ceci est
l'objet des propositions \ref{prop:z4} et \ref{prop:zk}.
La proposition \ref{prop:z0} calcule explicitement le terme $\Zp$
(modulo de nouveaux termes d'erreur qu'il s'agira de contr\^oler, {\it
  cf.} la sous-section \ref{subsec:zp} pour plus de d\'etails)
gr\^ace \`a l'expression exacte de  $\ecN_S(\bd ,\becE,\bfrE)$ pour $\bd$
\og grand \fg, et d\'egage ainsi le terme principal de la fonction z\^eta des hauteurs.

La d\'emonstration du th\'eor\`eme \ref{thm:resprinc} s'obtient alors en combinant les
propositions  \ref{prop:z4}, \ref{prop:zk} et \ref{prop:z0}.

\section{Quelques lemmes}

Nous rassemblons dans cette section quelques lemmes
qui nous seront utiles lors de la d\'emonstration du th\'eor\`eme
\ref{thm:resprinc}.

\subsection{Un lemme combinatoire}\label{subsec:combi}
\begin{nota}\label{nota:Fnu}
Soit $r\geq 1$ un entier.
Pour $\bnu\in \N^{\,r}$ posons
\begin{equation}
F_{\bnu}(\rho,\bT)
\eqdef
\sum_{\bn\in \N^{\,r}}\rho^{\Min(n_i+\nu_i)}\prod_i T_i^{\,n_i}\in
\Z[[\rho,(T_i)_{1\leq i\leq r}]]
\end{equation}
et
\begin{equation}
\wt{F}_{\bnu}(\rho,\bT)
\eqdef
(1-\rho\,\prod_i T_i)
\prod_i (1-T_i)\,F_{\bnu}(\rho,\bT).
\end{equation}
\end{nota}
Notons que pour tout $\mu\in\N$ on a
\begin{equation}\label{eq:form:mu}
F_{(\nu_i+\mu)}(\rho,\bT)=\rho^{\mu}\,F_{\bnu}(\rho,\bT).
\end{equation}

\begin{prop}\label{prop:coefFnu}
Soit $\bnu\in \N^{\,r}$.
\begin{enumerate}
\item\label{item:prop:coefFnu:1}
Soit $\bn\in \N^{\,r}$. 
On pose
\begin{align}
m&=\Min(n_i+\nu_i)\\
I_1&=\{i\in \{1,\dots,r\},\,\,n_i\geq 1\}\\
I_2&=\{i\in \{1,\dots,r\},\,\,n_i\geq 2\}\\
\text{et}\quad K&=\{i\in I_1,\,\,n_i+\nu_i\geq m+1\}.
\end{align}
Si $I_1\neq \{1,\dots,r\}$,
le coefficient d'indice $\bn$ de $\wt{F}_{\bnu}(\rho,\bT)$ 
vaut
\begin{equation}
\left\{
\begin{array}{ll}
\rho^{m}&\text{si}\quad I_1=\vide \\
\rho^m-\rho^{m-1}&\text{si}\quad I_1\neq \vide\quad \text{et}\quad
K=\vide \\
0&\text{si}\quad K\neq \vide .
\end{array}
\right.
\end{equation} 
Si $I_1=\{1,\dots,r\}$,
le coefficient d'indice $\bn$ de $\wt{F}_{\bnu}(\rho,\bT)$ 
vaut
\begin{equation}
\left\{
\begin{array}{ll}
0&\text{si}\quad I_2\cap K \neq \vide \\
-\rho^m&\text{si}\quad I_2=\vide\quad\text{et}\quad  K\neq \vide\\
0&\text{si}\quad I_2\neq\vide\quad\text{et}\quad  K=\vide\\
\rho^{m-1}-\rho^m &\text{si}\quad I_2\neq\vide\quad\text{et}\quad  K\neq
\vide\quad\text{et}\quad  I_2\cap K=\vide\\
-\rho^{m-1}&\text{si}\quad I_2=K=\vide.
\end{array}
\right.
\end{equation}
En particulier on a
\begin{equation}
\wt{F}_{(0,\dots,0)}(\rho,\bT)=1-\prod_i T_i.
\end{equation}
\item\label{item:prop:coefFnu:2}
$\wt{F}_{\bnu}(\rho,\bT)$ est un polyn\^ome dont le degr\'e partiel
en chaque $T_i$ est major\'e par $\Max(\nu_i)+1$.
\item\label{item:prop:coefFnu:3}
Soit $\rho\geq 1$,  $\eps>0$ et $(\eta_i)_{1\leq i\leq r}$ 
des nombres complexes de module $1$.
On a pour tout $\bnu$  la majoration
\begin{equation}
\abs{\wt{F}_{\bnu}\left(\rho,\left(\eta_i\,\rho^{-1}\right)\right)}
\leq (2+\Max(\nu_i)-\Min(\nu_i))^r\,\rho^{\,\Min(\nu_i)}
\end{equation}
\end{enumerate}
\end{prop}
\begin{demo}
Montrons le point \ref{item:prop:coefFnu:1}.
Notons $a_{\bn}$ le coefficient en question.

Supposons d'abord $I_1\neq \{1,\dots,r\}$. Un peu d'attention montre 
que $a_{\bn}$ s'\'ecrit alors
\begin{equation}
\sum_{0\leq k \leq r}(-1)^k
\sum_{
\substack{
J\subset I_1,\\
\card{J}=k
}
}
\rho^{\Min((n_i+\nu_i)_{i\notin J},(n_i-1+\nu_i)_{i\in
      J})}
\end{equation}
soit
\begin{equation}
a_{\bn}=
\sum_{0\leq k \leq r}(-1)^k
\sum_{
\substack{
J\subset I_1,\\
\card{J}=k
\\
J\setminus K\neq \vide
}
}
\rho^{m-1}
+
\sum_{0\leq k \leq r}(-1)^k
\sum_{
\substack{
J\subset K,\\
\card{J}=k
}
}
\rho^{m}
\end{equation}
On a donc, par un argument combinatoire classique,
\begin{equation}
a_{\bn}=\left\{
\begin{array}{ll}
\rho^m&\text{si}\quad I_1=\vide \\
\rho^m-\rho^{m-1}&\text{si}\quad I_1\neq \vide\quad \text{et}\quad
K=\vide \\
0&\text{si}\quad K\neq \vide.
\end{array}
\right.
\end{equation}

Supposons \`a pr\'esent $I_1=\{1,\dots,r\}$. 
Alors $a_{\bn}$ s'\'ecrit
\begin{align}
&\phantom{=}\sum_{0\leq k \leq r}(-1)^k
\sum_{
\substack{
J\subset \{1,\dots,r\},\\
\card{J}=k,
}
}
\rho^{\Min((n_i+\nu_i)_{i\notin J},(n_i-1+\nu_i)_{i\in
      J})}
\\
&\phantom{=}+
\sum_{0\leq k \leq r}(-1)^{k+1}
\sum_{
\substack{
J\subset I_2,\\
\card{J}=k,
}}
\rho^{1+\Min((n_i-1+\nu_i)_{i\notin J},(n_i-2+\nu_i)_{i\in
      J})}
\\
&=
\sum_{0\leq k \leq r}(-1)^k
\sum_{
\substack{
J\subset \{1,\dots,r\},\\
\card{J}=k,
\\
J\setminus I_2\neq \vide
}
}
\rho^{\Min((n_i+\nu_i)_{i\notin J},(n_i-1+\nu_i)_{i\in
      J})}.
\end{align}

On en d\'eduit qu'on a 
\begin{equation}\label{expr:an}
a_{\bn}=\sum_{0\leq k \leq r}(-1)^k
\sum_{
\substack{
J\subset K,\\
\card{J}=k \\
J \setminus I_2 \neq \vide
}
}
\rho^{m}
+
\sum_{0\leq k \leq r}(-1)^k
\sum_{
\substack{
J\subset I_1,\\
\card{J}=k
\\
J \setminus I_2 \neq \vide
\\
J\setminus  K \neq \vide
}
}
\rho^{m-1}.
\end{equation}
Or on a, toujours par un argument combinatoire classique,
\begin{equation}
\sum_{0\leq k \leq r}(-1)^k
\sum_{
\substack{
J\subset K,\\
\card{J}=k \\
J \setminus I_2 \neq \vide
}
}
\rho^{m}
=
\left\{
\begin{array}{ll}
0&\text{si}\,I_2\cap K \neq \vide \\
-\rho^m&\text{si}\,I_2=\vide\,\text{et}\, K\neq \vide\\
0&\text{si}\,I_2\neq\vide\,\text{et}\, K=\vide\\
-\rho^{m}&\text{si}\,I_2\neq\vide\,\text{et}\, K\neq
\vide\,\text{et}\, I_2\cap K=\vide\\
0&\text{si}\,I_2=K=\vide
\end{array}
\right.
\end{equation}
et 
\begin{equation}
\sum_{0\leq k \leq r}(-1)^k
\sum_{
\substack{
J\subset I_1,\\
\card{J}=k
\\
J \setminus I_2 \neq \vide
\\
J\setminus  K \neq \vide
}
}
\rho^{m-1}
=
\left\{
\begin{array}{ll}
0&\text{si}\,I_2\cap K \neq \vide \\
0&\text{si}\,I_2=\vide\,\text{et}\, K\neq \vide\\
0&\text{si}\,I_2\neq\vide\,\text{et}\, K=\vide\\
\rho^{m-1}&\text{si}\,I_2\neq\vide\,\text{et}\, K\neq
\vide\,\text{et}\, I_2\cap K=\vide\\
-\rho^{m-1}&\text{si}\,I_2=K=\vide
\end{array}
\right.
\end{equation}
On en d\'eduit le point \ref{item:prop:coefFnu:1}.

Montrons le point \ref{item:prop:coefFnu:2}.
Soit $\bn\in \N^d$ tel que le coefficient d'indice $\bn$
de $\wt{F}_{\bnu}(\rho,\bT)$ soit non nul. 

Supposons $I_1\neq \{1,\dots,r\}$.
Soit $i_0\in \{1,\dots,r\}\setminus I_1$.
On a donc
\begin{equation}
m\leq n_{i_0}+\nu_{i_0}=\nu_{i_0}
\end{equation}
Comme le coefficient d'indice $\bn$ est non nul, d'apr\`es le
point \ref{item:prop:coefFnu:1} on a 
\begin{equation}
\forall i\in I_1,\quad n_{i}+\nu_{i}=m\leq \nu_{i_0}
\end{equation}
On a donc 
\begin{equation}
\forall i\in \{1,\dots,r\},\quad n_i\leq \Max(\nu_j)
\end{equation}

Supposons \`a pr\'esent $I_1=\{1,\dots,r\}$.
Comme le coefficient d'indice $\bn$ est non nul, 
d'apr\`es le
point \ref{item:prop:coefFnu:1} on a 
 soit $I_2=\vide$ soit $I_2\neq \vide$
et $K\neq \vide$ et $K\cap I_2=\vide$.

Si $I_2=\vide$ alors $n_i=1$ pour tout $i$.

Supposons \`a pr\'esent 
$I_2\neq \vide$ et $K\neq \vide$ et $K\cap I_2=\vide$.
Comme $K\neq \vide$ et $K\cap I_2=\vide$, il existe $i_0\in I_1\setminus I_2$
tel que
\begin{equation}
1+\nu_{i_0}=n_{i_0}+\nu_{i_0}\geq m+1
\end{equation}
d'o\`u $m\leq \nu_{i_0}$.

Comme $K\cap I_2=\vide$, on a pour tout $i\in I_2$
\begin{equation}
n_i+\nu_i=m\leq \nu_{i_0}
\end{equation}
et donc
\begin{equation}
\forall i\in I_2,\quad n_i\leq \Max(\nu_j)
\end{equation}

Finalement on a montr\'e que quelle que soit la valeur de $\bn$
telle que le coefficient d'indice $\bn$ est non nul, on a
\begin{equation}
\forall i,\quad n_i\leq \Max(\nu_j)+1.
\end{equation}

Montrons le point \ref{item:prop:coefFnu:3}. 
En utilisant \eqref{eq:form:mu}, on se ram\`ene aussit\^ot 
au cas o\`u $\Min(\nu_i)=0$.
Dans ce cas, on a pour tout $\bn$
\begin{equation}
\Min(n_i+\nu_i)-\sum n_i\leq 0.
\end{equation}
D'apr\`es les point \ref{item:prop:coefFnu:1}
et \ref{item:prop:coefFnu:2}, 
on a 
\begin{align}
\abs{\wt{F}_{\bnu}(\rho,\eta_i\,\rho^{-1})}
&
\leq
\sum_{0\leq n_i\leq \Max(\nu_j)+1
}
2\,\rho^{\Min(n_i+\nu_i)-\sum_i n_i}
\\
&
\leq
2\, \Max(\nu_i+2)^r.
\end{align}
\end{demo}
\subsection{Une estimation}
\begin{lemme}\label{lm:estim:plusieurs:variables}
Soit $r\geq 1$ un entier, $(a_{\bn})\in \C^{\,\N^{\,r}}$ et $\rho>0$ un r\'eel.
On suppose qu'il existe un r\'eel $\eps>0$
tel que la s\'erie
\begin{equation}
F(\bz)
=
\prod_i (1-\,\rho\,z_i)
\sum_{\bn\in \N^{\,r}} a_{\bn}\,\prod_i z_i^{n_i}
\end{equation}
converge absolument dans le domaine $\abs{z_i}\leq \rho^{-1}+\eps$.

Soit
\begin{equation}
\norm{F}_{\rho^{-1}}
\eqdef
\Maxu{\abs{\eta_i}=1} \,\,\abs{F(\eta_i\,\rho^{-1})}.
\end{equation}
On a alors
\begin{equation}
\forall \bn\in \N^{\,r},\quad
\abs{a_{\bn}}\leq 
\prod_i (n_i+1)
\norm{F}_{\rho^{-1}}\,\rho^{\sumu{i} n_i}.
\end{equation}
\end{lemme}
\begin{demo}
Si on \'ecrit $F(\bz)=\sumu{\bn\in \N^{\,r}}b_{\bn}\,\produ{i} z_i^{n_i}$
on a d'apr\`es les estimations de Cauchy
\begin{equation}
\forall \bn\in \N^{\,r},\quad
\abs{b_{\bn}}\leq \norm{F}_{\rho^{-1}}\,\rho^{\sumu{i} n_i}.
\end{equation} 
Or on a 
\begin{equation}
a_{\bn}=\sum_{\bm+\bm'=\bn} b_{\bm}\,\rho^{\sumu{i} m'_i}.
\end{equation}
Le r\'esultat en d\'ecoule aussit\^ot.
\end{demo}

\begin{lemme}\label{lm:ctheta}
Soit $r\geq 1$ et $\becD\in \diveffc^r$. 
On pose pour $\bd\in\N^{\,r}$
\begin{equation}
a_{\bd,\becD}
\eqdef \sum_{\substack{
\becG\in \diveffc^r\\
\deg(\becG)=\bd
}}
q^{\,\deg(\pgcd(\ecD_i+\ecG_i))}.
\end{equation}

Il existe une constante $c>0$
(ne d\'ependant que de $\courbe$)
telle qu'on ait la propri\'et\'e suivante : pour tout r\'eel $\theta>0$,
il existe une constante $c_{\theta}>0$ (ne d\'ependant que de $q$) 
telle qu'on ait
\begin{equation}
\forall \bd\in \N^{\,r},
\quad 
\abs{a_{\bd,\becD}}
\leq 
c_{\theta}\,c^{1+r\,(2+\sumu{v} 2\,\Max(v(\ecD_i))-\Min(v(\ecD_i))}\,
q^{\,\deg(\pgcd(\ecD_i)}\,
q^{\,(1+\theta)\,\sumu{i} d_i}.
\end{equation}
et
\begin{equation}
\forall \bd\in \N^{\,r},
\quad 
\abs{a_{\bd,\becD}}
\leq 
c_{\theta}\,(d_1+1)\,c^{1+r\,(2+\sumu{v} 2\,\Max(v(\ecD_i))-\Min(v(\ecD_i))}\,
q^{\,\deg(\pgcd(\ecD_i)}\,
q^{\,d_1+(1+\theta)\,\sumu{2\leq i\leq r} d_i}.
\end{equation}
\end{lemme}
\begin{demo}
Formons la s\'erie g\'en\'eratrice
\begin{align}
Z_{\becD}(\bT)
&\eqdef\sum_{\bd\in \N^{\,r}}
a_{\bd,\becD}
\,\,\prod T_i^{d_i}\\
&=\sum_{\becG\in \diveffc^r}
q^{\,\deg(\pgcd(\ecD_i+\ecG_i))}
\,\,\prod T_i^{\deg(\ecG_i)}.
\end{align}
Cette s\'erie s'\'ecrit comme le produit eul\'erien
({\it cf.} les notations de la section \ref{subsec:combi})
\begin{equation}
\prod_{v\in \courbe^{(0)}}
F_{(v(\ecD_i))}(q_v,\bT^{f_v}).
\end{equation}

On a donc d'apr\`es la proposition \ref{prop:coefFnu}
\begin{equation}
Z_{\becD}(\bT)=\left(\prod_i
  Z_{\courbe}(T_i)\right) F(\bT) \,G_{\becD}(\bT)
\end{equation}
avec
\begin{equation}
F(\bT)=\prod_{v\in \courbe^{(0)}} 
\frac{1-(\produ{i} T_i)^{f_v}}{1-q_v\,(\produ{i} T_i)^{f_v}}
\end{equation}
et
\begin{equation}
G_{\becD}(\bT)=\prod_{\substack{
v\in \courbe^{(0)}\\
(v(\ecD_i))\neq (0)}}
\frac{\wt{F}_{(v(\ecD_i))}(q_v,\bT^{f_v})}{1-(\produ{i} T_i)^{f_v}}.
\end{equation}
Compte tenu du  point \ref{item:prop:coefFnu:3} de la proposition 
\ref{prop:coefFnu}, on a 
\begin{align}
\norm{G_{\becD}}_{q^{-1}}
&
\leq 
\prod_{
\substack{
v\in \courbe^{(0)}\\
(v(\ecD_i))\neq (0)}
}
\frac
{e^{r\,(2+\Max(v(\ecD_i))-\Min(v(\ecD_i)))}
q_v^{\Min(v(\ecD_i))}
}
{1-q^{-r}}
\\
&
\leq 
q^{r \card{
\{v\in \courbe^{(0)},
(v(\ecD_i))\neq (0)\}
}
}
\,
e^{r\,\sumu{v}(2+\Max(v(\ecD_i))-\Min(v(\ecD_i)))}
q^{\deg(\pgcd(\ecD_i))}
\\
&
\leq 
q^{r \sumu{v} \Max(v(\ecD_i))}
\,
e^{r\,\sumu{v}(2+\Max(v(\ecD_i))-\Min(v(\ecD_i)))}
q^{\deg(\pgcd(\ecD_i))}.
\end{align}
D'apr\`es le lemme \ref{lm:estim:plusieurs:variables}, 
on a le r\'esultat voulu.
\end{demo}

\subsection{Comptage de sections globales}\label{subsec:compt:sec}

Les r\'esultats de cette partie sont \`a la base de l'estimation
de la quantit\'e $\ecN_S(\bd,\becE,\bfrD)$ introduite \`a la
section \ref{subsec:prelim}.
Les d\'emonstrations reposent sur de l'alg\`ebre lin\'eaire \'el\'ementaire, 
ainsi que sur le th\'eor\`eme de Riemann-Roch (lemme \ref{lm:estim:dim}).

\begin{lemme}\label{lm:clef:part:1}
Soient $\ecD_1,\ecD_2,\ecD'_1$ et $\ecD'_2$ des diviseurs de $\courbe$
tels qu'on ait
\begin{equation}
\ecD_1+\ecD'_1\sim \ecD_2+\ecD'_2.
\end{equation}
Soit $s_1$ (respectivement $s_2$) une section globale non nulle de $\ecO_{\courbe}(\ecD_1)$
(respectivement $\ecO_{\courbe}(\ecD_2)$).
On fixe un isomorphisme 
\begin{equation}
\ecO_{\courbe}(\ecD_1+\ecD'_1)\isom
\ecO_{\courbe}(\ecD_2+\ecD'_2),
\end{equation}
ce qui permet
de d\'efinir l'application lin\'eaire
\begin{equation}
\varphi_{s_1,s_2}\,:\,H^0(\courbe,\ecO_{\courbe}(\ecD'_1))
\times 
H^0(\courbe,\ecO_{\courbe}(\ecD'_2))
\longto H^0(\courbe,\ecO_{\courbe}(\ecD_1+\ecD'_1))
\end{equation}
qui \`a $(t_1,t_2)$ associe $t_1\,s_1+t_2\,s_2$.

Soit $\delta$ le degr\'e de $\ecD_1+\ecD'_1$ (ou ce qui revient au m\^eme
celui de $\ecD_2+\ecD'_2$).
\begin{enumerate}
\item\label{item:lm:clef:part:1:1}
On suppose qu'on a l'in\'egalit\'e 
\begin{equation}
\delta<\deg(\ecD_1)+\deg(\ecD_2)-\deg(\pgcd[\ddiv(s_1),\ddiv(s_2)]).
\end{equation}
Alors $\varphi_{s_1,s_2}$ est injective.
\item\label{item:lm:clef:part:1:2}
On suppose qu'on a l'in\'egalit\'e 
\begin{equation}
\delta\geq \deg(\ecD_1)+\deg(\ecD_2)-\deg(\pgcd[\ddiv(s_1),\ddiv(s_2)]).
\end{equation}
Alors on a 
\begin{equation}
\dim(\Ker(\varphi_{s_1,s_2}))\leq 1+\delta-\deg(\ecD_1)-\deg(\ecD_2)+\deg(\pgcd[\ddiv(s_1),\ddiv(s_2)]).
\end{equation}
\item\label{item:lm:clef:part:1:3}
On suppose qu'on a l'in\'egalit\'e 
\begin{equation}
\delta\geq \deg(\ecD_1)+\deg(\ecD_2)-\deg(\pgcd(\ddiv(s_1),\ddiv(s_2)))+2\,g-1.
\end{equation}
Alors on a 
\begin{equation}
\dim(\Ker(\varphi_{s_1,s_2}))=1-g+\delta-\deg(\ecD_1)-\deg(\ecD_2)+\deg(\pgcd(\ddiv(s_1),\ddiv(s_2)))
\end{equation}
et 
\begin{equation}
\im(\varphi_{s_1,s_2})=\{s\in
H^0(\courbe,\ecO_{\courbe}(\ecD_1+\ecD'_1))\setminus\{0\},\,\ddiv(s)\geq
\pgcd(\ddiv(s_1),\ddiv(s_2))\}\cup \{0\}
\end{equation}
\end{enumerate}
\end{lemme}
\begin{demo}
Quitte \`a remplacer $\ecD'_2$ par un diviseur lin\'eairement \'equivalent,
on peut supposer qu'on a 
\begin{equation}
\ecD_1+\ecD'_1=\ecD_2+\ecD'_2\eqdef \ecD.
\end{equation}
On peut \'egalement supposer qu'on a $\pgcd(\ddiv(s_1),\ddiv(s_2))=0$.
Soit $(t_1,t_2)\in \Ker(\varphi_{s_1,s_2})$. Notons qu'on a $t_1=0$
si et seulement si $t_2=0$.
Supposons $(t_1,t_2)\neq (0,0)$. Comme on a
$s_1\,t_1=-s_2\,t_2$, on a
\begin{equation}
\div(t_1)+\div(s_1)=\div(t_2)+\div(s_2)
\end{equation}
d'o\`u 
\begin{equation}\label{eq:div}
\div(t_1)\geq \div(s_2).
\end{equation}
En particulier, on a 
\begin{equation}
\deg(\ecD_2)\leq(\ecD'_1)=\deg(\ecD)-\deg(\ecD_1) 
\end{equation}
ce qui montre la premi\`ere assertion.

Par ailleurs on en d\'eduit facilement de \eqref{eq:div} que l'application 
$u\mapsto (u\,s_2,u\,s_1)$ est un isomorphisme de
$H^0(\courbe,\ecO_{\courbe}(\ecD'_1-\ecD_2))$ sur
$\Ker(\varphi_{(s_1,s_2)})$.
On en d\'eduit les deux derni\`eres assertions, compte tenu du lemme
\ref{lm:estim:dim}.
\end{demo}

\begin{cor}\label{cor:clef:part}
Soient $\ecD_1,\ecD_2,\ecD_3,\ecD'_1,\ecD'_2$ et $\ecD'_3$ des diviseurs de $\courbe$
tels qu'on ait
\begin{equation}
\ecD_1+\ecD'_1\sim \ecD_2+\ecD'_2\sim \ecD_3+\ecD'_3.
\end{equation}
Soit $s_1$ (respectivement $s_2$, $s_3$) une section globale non nulle de $\ecO_{\courbe}(\ecD_1)$
(respectivement $\ecO_{\courbe}(\ecD_2)$, $\ecO_{\courbe}(\ecD_3)$).

On fixe des isomorphismes 
\begin{equation}
\ecO_{\courbe}(\ecD_1+\ecD'_1)\isom \ecO_{\courbe}(\ecD_2+\ecD'_2)
\end{equation}
et 
\begin{equation}
\ecO_{\courbe}(\ecD_1+\ecD'_1)\isom \ecO_{\courbe}(\ecD_3+\ecD'_3),
\end{equation} 
ce qui permet de d\'efinir l'application lin\'eaire
\begin{equation}
\varphi_{s_1,s_2,s_3}\,:\,H^0(\courbe,\ecO_{\courbe}(\ecD'_1))
\times H^0(\courbe,\ecO_{\courbe}(\ecD'_2))
\times H^0(\courbe,\ecO_{\courbe}(\ecD'_3))
\longto H^0(\courbe,\ecO_{\courbe}(\ecD_1+\ecD'_1))
\end{equation}
qui \`a $(t_1,t_2,t_3)$ associe $t_1\,s_1+t_2\,s_2+t_3\,s_3$.

On note $\delta$ le degr\'e de $\ecD_1+\ecD'_1$.
\begin{enumerate}
\item\label{item:cor:clef:part:2}
Le cardinal de l'ensemble des \'el\'ements $(t_1,t_2,t_3)$ de
$\Ker(\varphi_{s_1,s_2,s_3})$ v\'erifiant $t_3\neq 0$
est major\'e par 
\begin{equation}
q^{2+2\,\delta-\deg(\ecD_1)-\deg(\ecD_2)-\deg(\ecD_3)}+q^{1+\deg(\ecD'_3)}.
\end{equation}
\item\label{item:cor:clef:part:3}
On suppose qu'on a
\begin{equation}\label{eq:hyp:cor:clef:part:3}
\delta\geq \deg(\ecD_1)+\deg(\ecD_2)-1
\end{equation}
Le cardinal de l'ensemble des \'el\'ements $(t_1,t_2,t_3)$ de
$\Ker(\varphi_{s_1,s_2,s_3})$ v\'erifiant $t_3\neq 0$
est major\'e par 
\begin{equation}
q^{2+2\,\delta-\deg(\ecD_1)-\deg(\ecD_2)-\deg(\ecD_3)}.
\end{equation}
\item\label{item:cor:clef:part:1}
On suppose qu'on a
\begin{equation}\label{hyp:cor:clef:part}
\delta\geq \deg(\ecD_1)+\deg(\ecD_2)+2\,g-1
\end{equation}
et 
\begin{equation}\label{hyp:cor:clef:part:2}
\delta\geq \deg(\ecD_2)+\deg(\ecD_3)+2\,g-1.
\end{equation}
Alors le noyau de $\varphi_{s_1,s_2,s_3}$ est de dimension
\begin{equation}
2\,(\delta+1-g)-\sum_{i=1}^3 \deg(\ecD_i)+\deg(\pgcd(\ddiv(s_i)).
\end{equation}
\end{enumerate}
\end{cor}
\begin{demo}
On peut supposer qu'on a 
\begin{equation}
\ecD_1+\ecD'_1=\ecD_2+\ecD'_2=\ecD_3+\ecD'_3\eqdef \ecD
\end{equation}
et qu'on a 
$\pgcd(\ddiv(s_i))=0$.

Soit $t_3$
un \'el\'ement non nul 
de
$
H^0(\courbe,\ecO_{\courbe}(\ecD'_3))
$ 
qui est dans l'image 
de $\Ker(\varphi_{s_1,s_2,s_3})$
par la troisi\`eme projection.
On a donc 
\begin{equation}
\div(t_3)+\div(s_3)\geq \pgcd[\div(s_1),\div(s_2)]
\end{equation}
soit 
\begin{equation}
\div(t_3)\geq \pgcd[\div(s_1),\div(s_2)].
\end{equation}
D'apr\`es le lemme \ref{lm:estim:dim}, le cardinal de
l'ensemble des \'el\'ements v\'erifiant cette propri\'et\'e  
est major\'e par 
\begin{equation}
q^{1+\deg(\ecD'_3)-\deg(\pgcd[\div(s_1),\div(s_2)])}.
\end{equation}
Par ailleurs si $t_3$ est un tel \'el\'ement, l'ensemble 
des couples $(t_1,t_2)$ tels que $(t_1,t_2,t_3)\in
\Ker(\varphi_{s_1,s_2,s_3})$ est un espace affine de direction
$\Ker(\varphi_{s_1,s_2})$. D'apr\`es les points
\ref{item:lm:clef:part:1:1}
et \ref{item:lm:clef:part:1:2}
du lemme \ref{lm:clef:part:1}, le cardinal de cet ensemble est major\'e
par 
\begin{equation}
q^{1+\deg(\ecD)-\deg(\ecD_1)-\deg(\ecD_2)+\deg(\pgcd[\ddiv(s_1),\ddiv(s_2)])}
\end{equation}
si
\begin{equation}\label{eq:cond:posit}
1+\deg(\ecD)-\deg(\ecD_1)-\deg(\ecD_2)+\deg(\pgcd[\ddiv(s_1),\ddiv(s_2)])\geq
0
\end{equation}
(ce qui est toujours v\'erifi\'e si l'hypoth\`ese \ref{eq:hyp:cor:clef:part:3}
est satisfaite)
et par
\begin{equation}
q^{1+\deg(\ecD)-\deg(\ecD_1)-\deg(\ecD_2)+\deg(\pgcd[\ddiv(s_1),\ddiv(s_2)])}+1
\end{equation}
sinon. 
On en d\'eduit les points \ref{item:cor:clef:part:2}
et \ref{item:cor:clef:part:3}.

Montrons le point \ref{item:cor:clef:part:1}.
D'apr\`es le point \ref{item:lm:clef:part:1:3} du lemme \ref{lm:clef:part:1} et l'hypoth\`ese \eqref{hyp:cor:clef:part}
un \'el\'ement $t_3$ non nul
de
$
H^0(\courbe,\ecO_{\courbe}(\ecD'_3))
$ 
est dans l'image 
de $\Ker(\varphi_{s_1,s_2,s_3})$
par la troisi\`eme projection si et seulement si 
\begin{equation}
\div(t_3)+\div(s_3)\geq \pgcd[\div(s_1),\div(s_2)]
\end{equation}
{\it i.e.} si et seulement si 
\begin{equation}
\div(t_3)\geq \pgcd[\div(s_1),\div(s_2)].
\end{equation}
Ainsi l'image de la projection de $\Ker(\varphi_{s_1,s_2,s_3})$
sur $H^0(\ecO_{\courbe}(\ecD'_3))$ est isomorphe
\`a $H^0(\ecO_{\courbe}(\ecD'_3-\pgcd[\div(s_1),\div(s_2)]))$.

Or on a 
\begin{equation}
\deg(\ecD'_3-\pgcd[\div(s_1),\div(s_2)])\geq \deg(\ecD'_3)-\deg(\ecD_2)
\end{equation}
soit d'apr\`es l'hypoth\`ese \eqref{hyp:cor:clef:part:2}
\begin{equation}
\deg(\ecD'_3-\pgcd[\div(s_1),\div(s_2)])\geq 2\,g-1.
\end{equation}
On conclut gr\^ace au lemme \ref{lm:estim:dim}
et au point \ref{item:lm:clef:part:1:3} du lemme
\ref{lm:clef:part:1}.
\end{demo}

\begin{cor}\label{cor:nsO}
Soit $\bd\in \N^4$  et 
$\becE\in \diveffc^7$  v\'erifiant
\begin{equation}
\forall i\in \{1,2,3\},\quad
\psi_i(\bd,\becE )\geq 0.
\end{equation}
\begin{enumerate}
\item\label{item:cor:nsO}
Pour tout $\bfrE\in \piczcSe $,
$
\ecN_{S}(\bd,\becE,\bfrE)
$
est major\'e par 
\begin{multline}
q^{2+3\,d_0+\sumu{i}d_i+2\,\deg(\ecE_0)+\sumu{i}\deg(\ecE_i)-\sumu{i}\deg(\ecF_i)}
\sum_{(s_i)\in \produ{i} \Homogs_{\frE_i,d_i}}
q^{\deg(\pgcd(\ddiv(s_i)+\ecE_i+\ecF_i))}
\\
+q^{5+2\,d_0+2\,d_1+2\,d_2+d_3+\,\deg(\ecE_0)+\deg(\ecE_1)+\deg(\ecE_2)-\deg(\ecF_3)}
\end{multline}
\item 
On suppose qu'il existe $j\in \{1,2,3\}$
tel que la condition
\begin{equation}
\phi_j(\bd ,\becE )\geq 2\,g-1
\end{equation}
est v\'erifi\'ee. Alors
pour tout $\bfrE\in \piczcSe $,
$
\ecN_{S}(\bd,\becE,\bfrE)
$
est major\'e par 
\begin{equation}
q^{2+3\,d_0+\sumu{i}d_i+2\,\deg(\ecE_0)+\sumu{i}\deg(\ecE_i)-\sumu{i}\deg(\ecF_i)}
\sum_{(s_i)\in \produ{i} \Homogs_{\frE_i,d_i}}
q^{\deg(\pgcd(\ddiv(s_i)+\ecE_i+\ecF_i))}
\end{equation}
\item
On suppose qu'il existe $j,k\in \{1,2,3\}$ avec $j\neq k$ v\'erifiant
\begin{equation}
\phi_j(\bd ,\becE )\geq 2\,g-1
\end{equation}
et
\begin{equation}
\phi_k(\bd ,\becE )\geq 2\,g-1.
\end{equation}
Alors, pour tout $\bfrE\in \piczcSe$,
la quantit\'e $\ecN_{S,0}(\bd,\becE ,\bfrE)$ 
({\it cf.} la sous-section \ref{subsubsec:deczp})
est \'egale \`a 
\begin{multline}
\card{\Homogs_{\frE_0,d_0}}\,\,
q^{\,2(1-g)+2\,d_0+\sumu{i} d_i+2\,\deg(\ecE_0)+\sumu{i}\deg(\ecE_i)-\sumu{i}
  \deg(\ecF_i)}
\\
\times
\sum_{(s_i)\in \produ{i} \Homogs_{\frE_i,d_i}}
q^{\deg(\pgcd(\ddiv(s_i)+\ecE_i+\ecF_i))}
\end{multline}
\end{enumerate}
\end{cor}
\begin{demo}
Quitte \`a permuter les indices, on peut supposer qu'on a $j=2$ et $k=3$.
Pour tout $(s_i)\in \produ{i} \Homogs_{\frE_i,d_i}$
on applique alors les diff\'erents r\'esultats du  corollaire \ref{cor:clef:part}
avec 
\begin{equation}
\ecD_i=\frE_i+d_i\,\frD_1+\ecE_i+\ecF_i
\end{equation}
et
\begin{equation}
\ecD'_i=\frF_i
+
\left(
d_0+\sumu{j\neq i}d_i+\deg(\ecE_0)+
\sumu{j\neq i}\deg(\ecE_i)-\deg(\ecF_i)
\right)
\frD_1.
\end{equation}
On somme ensuite les contributions obtenues sur l'ensemble des
\'el\'ements $(s_i)\in \produ{i} \Homogs_{\frE_i,d_i}$.
Pour le point \ref{item:cor:nsO}, on utilise en outre le fait qu'on a,
d'apr\`es le lemme \ref{lm:estim:dim},
\begin{equation}
\card{\produ{1\leq i\leq 3} \Homogs_{\frE_i,d_i}}\leq q^{3+\sumu{i}d_i}.
\end{equation}
\end{demo}
En utilisant le lemme \ref{lm:ctheta}, on  d\'eduit aussit\^ot 
du corollaire \ref{cor:nsO}
le
corollaire suivant.
\begin{cor}\label{cor:majecns}
Il existe une constante $c>0$
(ne d\'ependant que de $\courbe$)
et pour tout r\'eel $\theta>0$,
une constante $c_{\theta}>0$ (ne d\'ependant que de $q$) telles qu'on ait
la propri\'et\'e suivante :
soit $\bd\in \N^4$ 
et $\becE\in \diveffc^7$ v\'erifiant
\begin{equation}
\forall i\in \{1,2,3\},\quad
\psi_i(\bd ,\becE )\geq 0.
\end{equation}
Alors
\begin{enumerate}
\item
pour tout $\bfrE\in \piczcSe $,
$
\ecN_{S}(\bd,\becE,\bfrE)
$
est major\'e par 
\begin{multline}
c_{\theta}\,c^{1+\sumu{v} 2\,\Max(v(\ecE_i)+v(\ecF_i))-\Min(v(\ecE_i)+v(\ecF_i))}
\\
\times\,q^{2+3\,d_0+(2+\theta)\,\sumu{i}d_i+2\,\deg(\ecE_0)+\sumu{i}\deg(\ecE_i)-\sumu{i}\deg(\ecF_i)+\deg(\pgcd(\ecE_i+\ecF_i))}
\\
+q^{5+2\,d_0+2\,d_1+2\,d_2+d_3+\deg(\ecE_0)+\deg(\ecE_1)+\deg(\ecE_2)-\deg(\ecF_3)};
\end{multline}
\item 
si en outre il existe $j\in \{1,2,3\}$ tel que la condition 
\begin{equation}
\phi_j(\bd ,\becE )\geq 2\,g-1
\end{equation}
est v\'erifi\'ee, pour tout $\bfrE\in \piczcSe$,
$
\ecN_{S}(\bd,\becE,\bfrE)
$
est major\'e par 
\begin{multline}
c_{\theta}\,c^{1+\sumu{v} 2\,\Max(v(\ecE_i)+v(\ecF_i))-\Min(v(\ecE_i)+v(\ecF_i))}
\\
\times\,(d_1+1)\,q^{2+3\,d_0+2\,d_1+(2+\theta)\,(d_2+d_3)+2\,\deg(\ecE_0)+\sumu{i}\deg(\ecE_i)-\sumu{i}\deg(\ecF_i)+\deg(\pgcd(\ecE_i+\ecF_i))}.
\end{multline}
\end{enumerate}
\end{cor}

\subsection{Quelques propri\'et\'es de la fonction $\mu_S$}
On note 
$\{0,1\}^7_S$
l'ensemble des \'el\'ements 
$(e_0,\be,\bbf)\in \{0,1\}^7$
v\'erifiant
\begin{equation}
\Min\left(
\left(\sumu{j\neq i}e_j+\sumu{j} f_j\right)_i,
\left(e_0+\sumu{j\neq i}(e_j+f_j)\right)_i,
e_0+\sumu{i}e_i
\right)=0.
\end{equation}

\begin{lemme}\label{lm:sum:un:imp:mus:zero}
\begin{enumerate}
\item
On a 
\begin{equation}
\forall\,(e_0,\be,\bbf)\in \{0,1\}^7,\quad
\ind_{\{0,1\}_S^7}(e_0,\be,\bbf)
=
\sum_{\substack{
0\leq e_0\leq e_0\\
0\leq e'_i\leq e_i\\
0\leq f'_i\leq f_i
}
}
\mu^0_S(e'_0,\be',\bbf')
\end{equation}
En particulier,
pour tout $(e_0,\be,\bbf)\in \{0,1\}^7_S\setminus (0,\dots,0)$,
$\mu^0_S(e_0,\be,\bbf)$ est nul.
\item
Soit $e_0\in \{0,1\}$. On a  
\begin{equation}
\sum_{(\be,\bbf)\in \{0,1\}^6}
\mu^0_S(e_0,\be,\bbf)
=0
\end{equation}
\item
Soit $(e_0,\be,\bbf)\in \{0,1\}^7$.
On suppose que l'une des conditions suivantes est v\'erifi\'ee
\begin{enumerate}
\item\label{item:1:lm:sum:un:imp:mus:zero}
$
e_0+\sumu{i} (e_i+f_i)=1$ ;
\item\label{item:2:lm:sum:un:imp:mus:zero}
$e_0=e_1=e_2=e_3=0$ ;
\item\label{item:3:lm:sum:un:imp:mus:zero}
$e_0=0$ et il existe un $i\in \{1,2,3\}$
tel que $e_i=1$, $f_i=0$ et $e_j=0$
pour $j\neq i$.
\end{enumerate}
Alors 
$
\mu^0_S(e_0,\be,\bbf)$ 
est nul.
\end{enumerate}
\end{lemme}
\begin{demo}
Le premier point n'est autre  
que le contenu de la remarque \ref{rq:mu0} dans le cas 
o\`u  $X=S$, compte tenu de la de la description de $\ecI_S$
donn\'ee \`a la section \ref{subsec:tors:p2}.
Les autres points s'en d\'eduisent aussit\^ot.
\end{demo}
\begin{lemme}\label{lm:prop:mus:3}
Pour $\eps>0$, la s\'erie 
\begin{multline}\label{eq:serie:3}
\sum_{\becE \in
\diveffc^{7}}
\abs{\mu_S(\becE) }\,q^{-\left(\frac{1}{2}+\eps\right)
\left[\deg(\ecE_0)+\sumu{i}\deg(\ecE_i)+\sumu{i}\deg(\ecF_i)\right]}
\end{multline}
est convergente.
\end{lemme}
\begin{demo}
D'apr\`es la proposition \ref{prop:mu},
la s\'erie en question s'\'ecrit comme le produit eul\'erien
\begin{equation}
\prod_{v\in \courbe^{(0)}} 
1+\sum_{
\substack{(e_0,\be,\bbf)\in \{0,1\}^7\\ (e_0,\be,\bbf)\neq (0,0,0)}}
\abs{\mu^0_S(e_0,\be,\bbf)}\,
q_v^{\,-(\frac{1}{2}+\eps)\left(e_0+\sumu{i} e_i+f_i\right)}
\end{equation}
Pour montrer que ce dernier produit converge, il suffit de montrer
que les conditions
\begin{equation}
(e_0,\be,\bbf)\neq (0,0,0)\text{ et }\mu^0_S(e_0,\be,\bbf)\neq 0
\end{equation}
entra\^\i nent 
\begin{equation}
-\left(\frac{1}{2}+\eps\right)\left(e_0+\sumu{i} e_i+f_i\right)<-1
\end{equation}
ce qui d\'ecoule aussit\^ot du point 3(a) du lemme \ref{lm:sum:un:imp:mus:zero}.
\end{demo}
\begin{lemme}\label{lm:prop:mus:1}
Les s\'eries
\begin{equation}
\sum_{\becE 
\in\diveffc^{7}
}
\abs{\mu_S(\becE) }\,
q^{-\frac{3}{2}\deg(\ecE_0)-\frac{3}{4}(\deg(\ecE_1)+\deg(\ecE_2))-\frac{3}{2}\deg(\ecE_3)+\frac{1}{8}(\deg(\ecF_1)+\deg(\ecF_2))
-\frac{3}{4}\deg(\ecF_3)}
\end{equation}
et
\begin{equation}
\sum_{\becE 
\in\diveffc^{7}
}
\abs{\mu_S(\becE) }\,
q^{-\frac{3}{2}\left(\deg(\ecE_0)+\sumu{i}\deg(\ecE_i)\right)}
\end{equation}
sont convergentes.
\end{lemme}
\begin{demo}
En raisonnant comme dans la preuve du lemme \ref{lm:prop:mus:3}, 
on voit qu'il suffit de montrer
que les conditions
\begin{equation}
(e_0,\be,\bbf)\neq (0,0,0)\text{ et }\mu^0_S(e_0,\be,\bbf)\neq 0
\end{equation}
entra\^\i nent 
\begin{equation}\label{eq:in:be0}
-\frac{3}{2}e_0-\frac{3}{4}(e_1+e_2)-\frac{3}{2}e_3
+\frac{1}{8}(f_1+f_2)
-\frac{3}{4}f_3<-1,
\end{equation}
respectivement
\begin{equation}\label{eq:in:be0:bis}
-\frac{3}{2}\,\left(e_0+\sum_i e_i\right)<-1.
\end{equation}
Ceci d\'ecoule facilement du point 3 du lemme \ref{lm:sum:un:imp:mus:zero}.

\end{demo}
\begin{lemme}\label{lm:prop:mus:2}
Soit $\rho>0$. La s\'erie 
\begin{multline}
\sum_{\becE \in
\diveffc^{7}}
\abs{\mu_S(\becE) }\,
\rho^{\sumu{v}
  2\,\Max(v(\ecE_i)+v(\ecF_i))-\Min(v(\ecE_i)+v(\ecF_i))}
\\\times
q^{\deg(\pgcd(\ecE_i+\ecF_i))}q^{-\frac{3}{4}\left[\deg(\ecE_0)+\sumu{i}\deg(\ecE_i)+\deg(\ecF_i)\right]}
\end{multline}
est convergente.
\end{lemme}
\begin{demo}
D'apr\`es la proposition \ref{prop:mu}, la s\'erie en question
s'\'ecrit comme le produit eul\'erien
\begin{equation}
\prod_{v\in \courbe^{(0)}} 
\left(1+
\!\!\!\!\!\!\sum_{
\substack{(e_0,\be,\bbf)\in \{0,1\}^7\\ (e_0,\be,\bbf)\neq (0,0,0)}}
\!\!\!\!\!\!
\abs{\mu^0_S(e_0,\be,\bbf)}\,\rho^{2\,\Min(e_i+f_i)-\Max(e_i+f_i)}\,
q_v^{-\frac{3}{4}\,\left[e_0+\sumu{i}(e_i+f_i)+\Min(e_i+f_i)\right]} 
\right)
\end{equation}
Pour montrer que ce dernier produit converge, il suffit de montrer
que les conditions
\begin{equation}
(e_0,\be,\bbf)\neq (0,0,0)\text{ et }\mu^0_S(e_0,\be,\bbf)\neq 0
\end{equation}
entra\^\i nent 
\begin{equation}
-\frac{3}{4}\,\left[e_0+\sumu{i}(e_i+f_i)\right]+\Min(e_i+f_i)<1
\end{equation}
ce qui, l\`a encore, d\'ecoule facilement du point 3 du 
lemme \ref{lm:sum:un:imp:mus:zero}.
\end{demo}

\section{D\'emonstration du r\'esultat principal}

\subsection{Le terme $Z_4$}
\begin{prop}\label{prop:z4}
Le rayon de convergence de la s\'erie
\begin{equation}
\sum_{\becE \in
\diveffc^{7}}
\mu_S(\becE)
\sum_{\bfrE\in \piczcSe }
Z_4(\becE ,\bfrE,T)
\end{equation}
est strictement sup\'erieur \`a $q^{-1}$.
\end{prop}
\begin{demo}
Rappelons que pour $\becE $ et 
$\bfrE$ donn\'es,
$Z_{4}(\becE ,\bfrE,T)$
est donn\'ee par l'expression
\begin{equation}
T^{3\,\deg(\ecE_0)+2\,\sumu{i}\deg(\ecE_i)}
\!\!\!\!\!\!\!\!
\sum_{
\substack{
d_0\geq 0\\
~
\\
d_i\geq 0,
\\
~
\\
\psi_i(\bd,\becE )\geq 0
\\
~
\\
\phi_{i}(\bd,\becE )<2\,g-1
}}
\!\!\!\!\!
\ecN_S\left(
\bd ,\becE,\bfrD\right)\,T^{3\,d_0+2\,\sumu{i}d_i}.
\end{equation}
Les conditions
\begin{equation}
\forall i\in \{1,2,3\},\quad \phi_{i}(\bd,\becE )<0
\end{equation}
entra\^\i nent les in\'egalit\'es
\begin{gather}
0\leq d_0<\deg(\ecF_1)+\deg(\ecF_2)+\deg(\ecF_3)+2\,g,
\\
0\leq d_1<\deg(\ecF_2)+\deg(\ecF_3)+2\,g,
\\
0\leq d_2<\deg(\ecF_1)+\deg(\ecF_3)+2\,g,
\\
0\leq d_3<\deg(\ecF_1)+\deg(\ecF_2)+2\,g.
\end{gather}
En particulier $Z_{4}(\becE ,\bfrE,T)$
est un polyn\^ome en $T$ \`a coefficients positifs.
La proposition \ref{prop:z4}
d\'ecoule alors du lemme \ref{lm:convz4}
ci-dessous.
\end{demo}
\begin{lemme}\label{lm:convz4}
Il existe un r\'eel $\theta>0$ tel que 
la s\'erie
\begin{equation}
\sum_{\becE \in
\diveffc^{7}}
\abs{\mu_S(\becE)}
\sum_{\bfrE\in \piczcSe }
Z_{4}(\becE ,\bfrE,q^{-1+\theta})
\end{equation}
soit convergente.
\end{lemme}
\begin{demo}
Fixons $\becE \in \diveffc^{7}$ et
$\bfrE\in \piczcSe $.
D'apr\`es le corollaire \ref{cor:majecns} et la remarque ci-dessus, 
on a
\begin{equation}\label{eq:maj:z4}
Z_{4}(\becE ,\bfrE,q^{-1+\theta})\leq A(\theta)+B(\theta)
\end{equation}
avec 
\begin{align}
A(\theta)
&
=c_{\theta}\,\,q^{2+(3\,\theta-1)\,\deg(\ecE_0)+(2\,\theta-1)\,\sumu{i}\deg(\ecE_i)-\sumu{i}\deg(\ecF_i)}
\notag
\\
&
\phantom{=}
\times q^{\deg(\pgcd(\ecE_i+\ecF_i))}\,c^{1+\sumu{v} 2\,\Max(v(\ecE_i)+v(\ecF_i))-\Min(v(\ecE_i)+v(\ecF_i))}
\notag
\\
&
\phantom{=}
\times \sum_{\substack{
0\leq d_0<\deg(\ecF_1)+\deg(\ecF_2)+\deg(\ecF_3)+2\,g
\\
0\leq d_1<\deg(\ecF_2)+\deg(\ecF_3)+2\,g
\\
0\leq d_2<\deg(\ecF_1)+\deg(\ecF_3)+2\,g
\\
0\leq d_3<\deg(\ecF_1)+\deg(\ecF_2)+2\,g
}}
q^{\,3\,\theta\,d_0+3\,\theta\,\sumu{i}d_i}
\end{align}
et
\begin{align}
B(\theta)
&=
q^{5+(3\,\theta-2)\,\deg(\ecE_0)+(2\,\theta-1)\,(\deg(\ecE_1)+\deg(\ecE_2))+(2\,\theta-2)\,\deg(\ecE_3)-\deg(\ecF_3)}
\notag
\\
&
\phantom{=}
\times
\sum_{\substack{
0\leq d_0<\deg(\ecF_1)+\deg(\ecF_2)+\deg(\ecF_3)+2\,g
\\
0\leq d_1<\deg(\ecF_2)+\deg(\ecF_3)+2\,g
\\
0\leq d_2<\deg(\ecF_1)+\deg(\ecF_3)+2\,g
\\
0\leq d_3<\deg(\ecF_1)+\deg(\ecF_2)+2\,g
}}
q^{\,3\,\theta\,d_0+2\,\theta\,d_1+2\,\theta\,d_2+(-1+2\,\theta)\,d_3}
\end{align}

En utilisant, pour $\rho>1$ et $N\geq 1$, la majoration $\sumu{0\leq d
<N}\rho^d\leq \frac{\rho^{N}}{\rho-1}$, on en d\'eduit les majorations
\begin{align}
A(\theta)
&
\leq c'_{\theta}\,q^{\,2+24\,g\,\theta+(3\,\theta-1)\,\deg(\ecE_0)
+(2\,\theta-1)\,\sumu{i}\deg(\ecE_i)+(9\,\theta-1)\sumu{i}\deg(\ecF_i)}
\notag
\\
&
\phantom{leq}
\times q^{\deg(\pgcd(\ecE_i+\ecF_i))}\,c^{1+\sumu{v} 2\,\Max(v(\ecE_i)+v(\ecF_i))-\Min(v(\ecE_i)+v(\ecF_i))}.
\label{eq:majatheta}
\end{align}
et

\begin{align}
B(\theta)&
\leq
c''_{\theta}\,q^{5+18\,g\,\theta+(3\,\theta-2)\,\deg(\ecE_0)+(2\,\theta-1)\,\left[\deg(\ecE_1)+\deg(\ecE_2)\right]+(2\,\theta-2)\,\deg(\ecE_3)}
\notag
\\
&
\phantom{leq}
\times q^{7\,\theta\,\deg(\ecF_1)+7\,\theta\,\deg(\ecF_2)+(-1+7\,\theta)\,\deg(\ecF_3)}.
\label{eq:majbtheta}
\end{align}
o\`u $c_{\theta}$ et $c''_{\theta}$ sont des constantes ne d\'ependant que
de $\theta$.
Les majorations \eqref{eq:maj:z4}, \eqref{eq:majatheta} et
\eqref{eq:majbtheta} ainsi que les  lemmes \ref{lm:prop:mus:1} et \ref{lm:prop:mus:2} 
montrent le lemme \ref{lm:convz4}.
\end{demo}

\subsection{Les termes $Z_{k}$ pour $1\leq k\leq 3$}
\begin{prop}\label{prop:zk}
Pour $1\leq k \leq 3$ la s\'erie 
\begin{equation}
\sum_{\becE \in
\diveffc^{7}}
\mu_S(\becE) 
\sum_{\bfrE\in \piczcSe }
Z_{k}(\becE ,\bfrE,T)
\end{equation}
est $q^{-1}$-contr\^ol\'ee \`a l'ordre $2$.
\end{prop}
\begin{demo}
Nous traitons le cas o\`u $k=3$, les autres cas s'en d\'eduisent par
permutation des indices.

Rappelons que
pour $\becE \in
\diveffc^{7}$ et 
$\bfrE\in \piczcSe $,
$
Z_{3}(\becE ,\bfrE,T)
$
est donn\'ee par l'expression
\begin{equation}
T^{3\,\deg(\ecE_0)+2\,\sumu{i}\deg(\ecE_i)}
\sum_{
\substack{
\bd\in \N^4
\\
~
\\
\forall 1\leq i\leq 3,\quad \psi_i(\bd,\becE )\geq 0
\\
~
\\
\phi_{1}(\bd,\becE )< 2\,g-1
\\
\phi_{2}(\bd,\becE )< 2\,g-1
\\
\phi_{3}(\bd,\becE )\geq  2\,g-1.
}}
\!\!\!\!\!
\ecN_S\left(\bd,\becE,\bfrD\right)
\,T^{3\,d_0+2\,\sumu{i}d_i}
\end{equation}
Soit $\theta$ un r\'eel strictement positif.
D'apr\`es le corollaire \ref{cor:majecns}, 
$Z_{3}(\becE ,\bfrE,T)$ est major\'ee par la s\'erie
\begin{align}
&\phantom{\times}c_{\theta}\,c^{1+\sumu{v} 2\,\Max(v(E_i)+v(F_i))-\Min(v(E_i)+v(F_i))}
\\
&\times T^{3\,\deg(\ecE_0)+2\,\sumu{i}\deg(\ecE_i)}\,
q^{2+2\,\deg(\ecE_0)+\sumu{i}\deg(\ecE_i)-\sumu{i}\deg(\ecF_i)+\deg(\pgcd(\ecE_i+\ecF_i))}
\\
&\times
\sum_{
\substack{
d_3\geq 0
\\
0\leq d_0<\frac{1}{2}\deg(\ecF_1)+\frac{1}{2}\deg(\ecF_2)+\deg(\ecF_3)+2\,g
\\
0\leq d_1<\deg(\ecF_2)+\deg(\ecF_3)+2\,g
\\
0\leq d_3<\deg(\ecF_3)+\deg(\ecF_3)+2\,g
}
}
(d_3+1)\,q^{3\,d_0+(2+\theta)\,(d_1+d_2)}\,T^{3\,d_0+2\,\sumu{i}d_i}
\end{align}
On en d\'eduit que $Z_{3}(\becE ,\bfrE,T)$
est major\'ee par la s\'erie 
\begin{equation}
\left(\sum_{d_3\geq 0} (d_3+1)(q\,T)^{d_3}\right)
\,
\wt{Z}_{3,\theta}(\becE ,T)
=
\frac{1}{(1-q\,T)^2}
\wt{Z}_{3,\theta}(\becE ,T)
\end{equation}
o\`u $\wt{Z}_{3,\theta}(\becE ,T)$
est donn\'ee par l'expression
\begin{align}
&c_{\theta}\,c^{1+\sumu{v} 2\,\Max(v(\ecE_i)+v(\ecF_i))-\Min(v(\ecE_i)+v(\ecF_i))}
\\
&T^{3\,\deg(\ecE_0)+2\,\sumu{i}\deg(\ecE_i)}\,
q^{2\,\deg(\ecE_0)+\sumu{i}\deg(\ecE_i)-\sumu{i}\deg(\ecF_i)+\deg(\pgcd(\ecE_i+\ecF_i))}
\\
&\sum_{
\substack{
0\leq d_0<\frac{1}{2}\,\deg(\ecF_1)+\frac{1}{2}\,\deg(\ecF_2)+\deg(\ecF_3)+2\,g
\\
0\leq d_1<\deg(\ecF_2)+\deg(\ecF_3)+2\,g
\\
0\leq d_2<\deg(\ecF_1)+\deg(\ecF_3)+2\,g
}
}
q^{3\,d_0+(2+\theta)\,(d_1+d_2)}\,T^{3\,d_0+2\,d_1+2\,d_2}
\end{align}
En raisonnant comme dans la preuve du lemme \ref{lm:convz4}, 
on montre
qu'il existe un r\'eel $\theta>0$ tel que 
la s\'erie
\begin{equation}
\sum_{\becE \in
\diveffc^{7}}
\abs{\mu_S(\becE)}
\wt{Z}_{3,\theta}(\becE ,q^{-1+\theta})
\end{equation}
soit convergente, ce qui conclut la preuve de la proposition \ref{prop:zk}.
\end{demo}

\subsection{Le terme $\Zp$}\label{subsec:zp}

\begin{prop}\label{prop:z0}
Il existe une s\'erie $\wt{\Zp}(T)$ de rayon de convergence strictement
sup\'erieur \`a $q^{-1}$ telle que 
\begin{equation}
\wt{\Zp}\left(q^{-1}\right)=
(q-1)^4\,q^{\,2(1-g)}
\prod_{v\in\courbe^{(0)}}
(1-q_v^{-1})^{\rg(\Pic(S))}\,\frac{\card{S(\kappa_{v})}}{q_v^{\,\dim(S)}}
\end{equation} 
et
\begin{equation}
\Zp(T)-Z_{\courbe}(q^2\,T^3)\,Z_{\courbe}(q\,T^2)^3\,\wt{\Zp}(T)
\end{equation}
est $q^{-1}$-contr\^ol\'ee \`a l'ordre $2$.
\end{prop}
La d\'emonstration de cette proposition occupe le reste de cette section.
\subsubsection{D\'ecomposition de $\Zp$}\label{subsubsec:deczp}
Pour $\becE \in \diveffc^7$
et $\bfrE\in \piczcSe $,
rappelons que $\Zp(\becE ,\bfrE,T)$
est donn\'ee par l'expression
\begin{equation}\label{eq:def:zp}
T^{3\,\deg(\ecE_0)+2\,\sumu{i}\deg(\ecE_i)}
\sum_{
\substack{
\bd\in \N^4
\\
~
\\
\forall 1\leq i\leq 3,\quad \psi_i(\bd ,\becE )\geq 0
\\
~
\\
\exists k,k',\,k\neq k',\,\phi_k(\bd ,\becE )\geq 0
\\
\text{et }\phi_{k'}(\bd ,\becE )\geq 0
}}
\!\!\!\!\!\!\!
\ecN_S\left(\bd, \becE, \bfrE\right)
\,T^{3\,d_0+2\,\sumu{i}d_i}.
\end{equation}

Soit $\Zpr{0}(\becE ,\bfrE,T)$ la s\'erie d\'efinie par l'expression
\eqref{eq:def:zp} o\`u l'on a remplac\'e
$\ecN_S$ par $\ecN_{S,0}$, o\`u $\ecN_{S,0}(\bd,\becE ,\bfrE)$ repr\'esente le
cardinal de l'ensemble des \'el\'ements 
\begin{equation}
(s_0,s_i,t_i)\in 
\Homogs_{\frE_0,d_0} 
\times
\prod_i \Homogs_{\frE_i,d_i} 
\times
\prod_i \Homog_{\frF_i,d_0+\sumu{j\neq i}d_j+\deg(\ecE_0)+\sumu{j\neq i}\deg(\ecE_j)-\deg(\ecF_i)}
\end{equation}
v\'erifiant la relation $\sumu{i} s_i\,t_i\,\scan{\ecE_i}\,\scan{\ecF_i}=0$.

Pour $1\leq k\leq 3$, soit $\Zpr{k}(\becE ,\bfrE,T)$ la s\'erie d\'efinie par l'expression
\eqref{eq:def:zp} 
o\`u l'on a remplac\'e
$\ecN_S$ par $\ecN_{S,k}$, o\`u $\ecN_{S,k}(\bd,\becE,\bfrE)$ repr\'esente le
cardinal de l'ensemble des \'el\'ements 
\begin{equation}
(s_0,(s_i),(t_i)_{i\neq k})\in 
\Homogs_{\frE_0,d_0} 
\times
\prod_i \Homogs_{\frE_i,d_i} 
\times
\prod_{i\neq k} \Homog_{\frF_i,d_0+\sumu{j\neq i}d_j+\deg(\ecE_0)+\sumu{j\neq i}\deg(\ecE_j)-\deg(\ecF_i)}
\end{equation}
v\'erifiant la relation $\sumu{i\neq k} s_i\,t_i\,\scan{\ecE_i}\,\scan{\ecF_i}=0$.

Soit enfin $\Zpr{4}(\becE ,\bfrE,T)$ la s\'erie d\'efinie par l'expression
\eqref{eq:def:zp} 
o\`u l'on a remplac\'e
$\ecN_S$ par $\ecN_{S,4}$, o\`u 
\begin{equation}
\ecN_{S,4}(\bd,\becE ,\bfrE)\eqdef\card{\Homogs_{\frE_0,d_0}}\,\produ{i}\card{\Homogs_{\frE_i,d_i}}.
\end{equation}
On a ainsi
\begin{equation}\label{eq:dec:nso}
\ecN_{S,0}(\bd,\becE ,\bfrE)=\ecN_{S}(\bd,\becE ,\bfrE)+\sum_{k=1}^3
\ecN_{S,k}(\bd,\becE ,\bfrE)-2\,\ecN_{S,4}(\bd,\becE ,\bfrE)
\end{equation}
d'o\`u l'\'ecriture 
\begin{equation}
\Zp(\becE ,T)
=
\Zpr{0}(\becE ,\bfrE,T)
-\sum_{k=1}^3 \Zpr{k}(\becE ,\bfrE,T)
+2\,\Zpr{4}(\becE ,\bfrE,T).
\end{equation}
La proposition \ref{prop:z0} d\'ecoule alors des propositions
\ref{prop:z04}, \ref{prop:z0k}, \ref{prop:z0err}
et \ref{prop:z0princ}.

\begin{rem}
Compte tenu de la description du
torseur universel de $S$ donn\'e \`a la section \ref{subsec:tors:p2}, 
et rappelant que $S_0$ d\'esigne la surface $S$ priv\'ee des droites
$(\ecE_i)_{0\leq i\leq 3}$ et $(\ecF_i)_{1\leq i\leq 3}$,
on voit
que la d\'ecomposition \eqref{eq:dec:nso} correspond 
\`a la d\'ecomposition
g\'eom\'etrique
\begin{equation}
S\setminus \cupu{0\leq i\leq 3} \ecE_i
=
S_0
\,
\sqcup
\,
\cupu{1\leq i\leq 3} \big(\ecF_i\setminus  \ecE_i\big).
\end{equation} 
\end{rem}

\subsubsection{Le terme $\Zpr{4}$}

\begin{prop}\label{prop:z04}
Le rayon de convergence de la s\'erie 
\begin{equation}
\sum_{\becE \in
\diveffc^{7}}
\mu_S(\becE)
\sum_{\bfrE\in \piczcSe }
\Zpr{4}(\becE ,\bfrE, T)
\end{equation}
est strictement sup\'erieur \`a $q^{-1}$.
\end{prop}
\begin{demo}
Rappelons que $\Zpr{4}(\becE ,\bfrE,T)$
est donn\'ee par l'expression
\begin{multline}
T^{3\,\deg(\ecE_0)+2\,\sumu{i}\deg(\ecE_i)}
\sum_{
\substack{
\bd\in \N^4\\
\\
\forall 1\leq i\leq 3,\quad \psi_i(\bd ,\becE )\geq 0
\\
~
\\
\exists k,k',\,k\neq k',\,\phi_k(\bd ,\becE )\geq 0
\\
\text{et }\phi_{k'}(\bd ,\becE )\geq 0
}}
\!\!\!\!\!\!\!\!\!\!\!\!
\card{\Homogs_{\frE_0,d_0}}\,\produ{i}\card{\Homogs_{\frE_i,d_i}}
\,T^{3\,d_0+2\,\sumu{i}d_i}
\end{multline}
D'apr\`es le lemme \ref{lm:estim:dim}, pour tout r\'eel $\theta$  strictement positif, 
$\Zpr{4}(\becE ,q^{-1+\theta})$
est major\'ee par 
\begin{equation}
q^{4+(3\theta-3)\deg(\ecE_0)+(2\,\theta-2)\,\sumu{i}\deg(\ecE_i)}
\sum_{\bd\in \N^4}q^{(3\,\theta-2)\,d_0+(2\,\theta-1)\,\sumu{i}d_i}.
\end{equation}
Ainsi pour $\theta$ assez petit on a 
\begin{equation}
\Zpr{4}(\becE ,\bfrE,q^{-1+\theta})
\leq
\frac{1}{(1-q^{-\frac{1}{2}})^4}
\,\,
q^{4+(3\theta-3)\deg(\ecE_0)+(2\,\theta-2)\,\sumu{i}\deg(\ecE_i)}.
\end{equation}
Le lemme \ref{lm:prop:mus:1} permet de conclure.
\end{demo}
\subsubsection{Les termes $\Zpr{k}$ pour $1\leq k\leq 3$}
\begin{prop}\label{prop:z0k}
Pour $1\leq k \leq 3$ la s\'erie 
\begin{equation}
\sum_{\becE \in
\diveffc^{7}}
\mu_S(\becE)
\sum_{\bfrE\in \piczcSe }
\Zpr{k}(\becE ,\bfrE,T)
\end{equation}
est $q^{-1}$-contr\^ol\'ee \`a l'ordre 2.
\end{prop}
\begin{demo}
Nous traitons le cas o\`u $k=1$, les autres cas s'en d\'eduisent par
permutation des variables.

Rappelons que $\Zpr{1}(\becE ,\bfrE, T)$
est donn\'ee par l'expression
\begin{equation}
T^{3\,\deg(\ecE_0)+2\,\sumu{i}\deg(\ecE_i)}
\sum_{
\substack{
\bd\in \N^4
\\
~
\\
\forall 1\leq i\leq 3,\quad \psi_i(\bd ,\becE )\geq 0
\\
~
\\
\exists k,k',\,k\neq k',\,\phi_k(\bd ,\becE )\geq 0
\\
\text{et }\phi_{k'}(\bd ,\becE )\geq 0
}}
\!\!\!\!\!\!\!\!\!\!\!\!
\ecN_{S,1}\left(
\bd,\becE,\bfrE
\right)\,
T^{3\,d_0+2\,\sumu{i}d_i}
\end{equation}
D'apr\`es le lemme \ref{lm:majnS1} ci-dessous, 
$\Zpr{1}(\becE ,\bfrE, T)$ est major\'ee par 
la s\'erie
\begin{multline}
T^{3\,\deg(\ecE_0)+2\,\sumu{i}\deg(\ecE_i)}
\sum_{
\bd\in \N^4
}
q^{5+ 2\,d_0+2\,d_1+2\,d_2+d_3+\deg(\ecE_0)+\deg(\ecE_1)+\deg(\ecE_2)-\deg(\ecF_3)}\,T^{3\,d_0+2\,\sumu{i}d_i},
\end{multline}
donc par la s\'erie
\begin{equation}
\frac{1}{(1-q\,T)^2}\wt{\Zpr{1}}(\becE, T),
\end{equation}
o\`u $\wt{\Zpr{1}}(\becE, T)$ est donn\'ee par l'expression
\begin{multline}
T^{3\,\deg(\ecE_0)+2\,\sumu{i}\deg(\ecE_i)}
\sum_{
d_0,d_3\geq 0}
q^{2\,d_0+d_3+\deg(\ecE_0)+\deg(\ecE_1)+\deg(\ecE_2)-\deg(\ecF_3)}\,T^{3\,d_0+2\,d_3}.
\end{multline}
Pour $\theta>0$ assez petit,  
$\wt{\Zpr{1}}(\becE ,q^{-1+\theta})$ est  major\'e
par 
\begin{equation}
\frac{1}
{(1-q^{-1})^2}
\,\,
q^{(3\theta-2)\,\deg(\ecE_0)+(2\theta-1)\deg(\ecE_1)+(2\theta-1)\deg(\ecE_2)+(2\theta-2)\,\deg(\ecE_3)-\deg(\ecF_3)}
.
\end{equation}
Le lemme \ref{lm:prop:mus:1} montre alors que
le rayon de convergence de la s\'erie 
\begin{equation}
\sum_{\becE \in
\diveffc^{7}}
\abs{\mu_S(\becE)}
\wt{\Zpr{1}}(\becE ,T)
\end{equation}
est strictement sup\'erieur \`a $q^{-1}$, d'o\`u le r\'esultat.
\end{demo}

\begin{lemme}\label{lm:majnS1}
Soit $\bd\in \N^4$  et 
$\becE\in \diveffc^7$  v\'erifiant
\begin{equation}
\forall i\in \{1,2,3\},\quad
\psi_i(\bd,\becE )\geq 0.
\end{equation}
On a pour tout $\bfrE\in \piczcSe$ la majoration
\begin{equation}
\ecN_{S,1}\left(\bd,\becE,\bfrE \right)
\leq
q^{\,5+2\,d_0+2\,d_1+2\,d_2+d_3+\deg(\ecE_0)+\deg(\ecE_1)+\deg(\ecE_2)-\deg(\ecF_3)}.
\end{equation}
\end{lemme}
\begin{demo}
Rappelons que $\ecN_{S,1}(\bd,\becE,\bfrE)$ d\'esigne 
le
cardinal de l'ensemble des \'el\'ements 
\begin{equation}
(s_0,(s_i),(t_2,t_3))\in 
\Homogs_{\frE_0,d_0} 
\times
\prod_{1\leq i\leq 3} \Homogs_{\frE_i,d_i} 
\times
\prod_{2\leq i\leq 3} \Homog_{\frF_i,d_0+\sumu{j\neq i}d_j+\deg(\ecE_0)+\sumu{j\neq i}\deg(\ecE_j)-\deg(\ecF_i)}
\end{equation}
v\'erifiant la relation 
\begin{equation}\label{eq:rel}
s_2\,t_2\,\scan{\ecE_2}\,\scan{\ecF_2}=s_3\,t_3\,\scan{\ecE_3}\,\scan{\ecF_3}.
\end{equation}
Si on fixe $(s_0,(s_i),t_3)$, il existe au plus un \'el\'ement $t_2$
satisfaisant la relation \eqref{eq:rel}. 
Ainsi
$\ecN_{S,1}(\bd,\becE,\bfrE)$
est major\'e par 
\begin{equation}
\card{\Homogs_{\frE_0,d_0}}
\prod_{1\leq i\leq 3} \card{\Homogs_{\frE_i,d_i}}
\,\card{\Homog_{\frF_3,d_0+d_1+d_2+\deg(\ecE_0)+\deg(\ecE_1)+\deg(\ecE_2)-\deg(\ecF_3)}}
\end{equation}
Le lemme \ref{lm:estim:dim} permet de conclure.
\end{demo}

\subsubsection{Le terme $\Zpr{0}$}
Rappelons que $\Zpr{0}(\becE ,\bfrE, T)$ est donn\'ee par l'expression
\begin{equation}
T^{3\,\deg(\ecE_0)+2\,\sumu{i}\deg(\ecE_i)}
\sum_{
\substack{
\bd\in \N^4\\
~
\\
\forall 1\leq i\leq 3,\quad \psi_i(\bd ,\becE )\geq 0
\\
~
\\
\exists k,k',\,k\neq k',\,\phi_k(\bd ,\becE )\geq 0
\\
\text{et }\phi_{k'}(\bd ,\becE )\geq 0
}}
\ecN_{S,0}\left(\bd,\becE,\bfrE\right)
\,
T^{3\,d_0+2\,\sumu{i}d_i}.
\end{equation}
Ainsi, d'apr\`es le corollaire \ref{cor:nsO}, 
$\Zpr{0}(\becE ,\bfrE, T)$
peut s'\'ecrire 
\begin{multline}
T^{3\,\deg(\ecE_0)+2\,\sumu{i}\deg(\ecE_i)}
\,
q^{2\,\deg(\ecE_0)+\sumu{i}\deg(\ecE_i)-\sumu{i} \deg(\ecF_i)}
\\ 
\times
\sum_{
\substack{
\bd\in \N^4
\\
~
\\
\forall 1\leq i\leq 3,\quad \psi_i(\bd ,\becE )\geq 0
\\
~
\\
\exists k,k',\,k\neq k',\,\phi_k(\bd ,\becE )\geq 0
\\
\text{et }\phi_{k'}(\bd ,\becE )\geq 0
}}
\!\!\!\!\!\!\!\!\!\!\!\!
\card{\Homogs_{\frE_0,d_0}}
\,\,
q^{2(1-g)+2\,d_0+\sumu{i} d_i}
\sum_{(s_i)\in \produ{i} \Homogs_{\frE_i,d_i}}
\!\!\!\!q^{\deg(\pgcd(\ddiv(s_i)+\ecE_i+\ecF_i))}
\,T^{3\,d_0+2\,\sumu{i}d_i}.
\end{multline}
On d\'efinit $\Zpr{0,\text{princ}}(\becE ,\bfrE, T)$
par l'expression 
\begin{multline}
T^{3\,\deg(\ecE_0)+2\,\sumu{i}\deg(\ecE_i)}
\,
q^{2\,\deg(\ecE_0)+\sumu{i}\deg(\ecE_i)-\sumu{i} \deg(\ecF_i)}
\\
\times 
\sum_{
\bd\in \N^4}
\card{\Homogs_{\frE_0,d_0}}
\,\,
q^{2(1-g)+2\,d_0+\sumu{i} d_i}
\sum_{(s_i)\in \produ{i} \Homogs_{\frE_i,d_i}}
\!\!\!\!q^{\deg(\pgcd(\ddiv(s_i)+\ecE_i+\ecF_i))}
\,T^{3\,d_0+2\,\sumu{i}d_i}.
\end{multline}
et on pose 
\begin{equation}
\Zpr{0,\text{err}}\eqdef\Zpr{0,\text{princ}}-\Zpr{0}.
\end{equation}

\subsubsection{Le terme $\Zpr{0,\text{err}}$}
\begin{prop}\label{prop:z0err}
La s\'erie 
\begin{equation}
\sum_{\becE \in
\diveffc^{7}}
\mu_S(\becE) \,
\sum_{\bfrE\in \piczcSe }
\Zpr{0,\text{err}}(\becE ,\bfrE, T)
\end{equation}
est $q^{-1}$-contr\^ol\'ee \`a l'ordre $2$.
\end{prop}
\begin{demo}
Par d\'efinition,
$\Zpr{0,\text{err}}(\becE ,\bfrE, T)$
est major\'ee par la somme des six s\'eries
\begin{multline}
T^{3\,\deg(\ecE_0)+2\,\sumu{i}\deg(\ecE_i)}
\,
q^{2\,\deg(\ecE_0)+\sumu{i}\deg(\ecE_i)-\sumu{i} \deg(\ecF_i)}
\\ 
\times
\sum_{
\substack{
\bd \in \ecA
}}
\card{\Homogs_{\frE_0,d_0}}
\,\,
q^{2(1-g)+2\,d_0+\sumu{i} d_i}
\sum_{(s_i)\in \produ{i} \Homogs_{\frE_i,d_i}}
\!\!\!\!q^{\deg(\pgcd(\ddiv(s_i)+\ecE_i+\ecF_i))}
\,T^{3\,d_0+2\,\sumu{i}d_i},
\end{multline}
o\`u $\ecA$ est le sous-ensemble de $\N^4$ constitu\'e des \'el\'ements
v\'erifiant successivement l'une des six conditions suivantes :
\begin{equation}
\psi_i(\bd,\becE)<0,\quad 1\leq i \leq 3,
\end{equation}
\begin{equation}
\phi_i(\bd,\becE)<0,\quad 1\leq i \leq 3.
\end{equation}
Nous nous contentons de d\'emontrer que 
la s\'erie $\sum_{\becE}
\mu_S(\becE) \,
\sum_{\bfrE}
\Zpr{0,\text{err},\psi_1}(\becE ,\bfrE, T)$,
o\`u
$\Zpr{0,\text{err},\psi_1}(\becE ,\bfrE, T)$ est 
donn\'ee par l'expression
\begin{multline}
T^{\,3\,\deg(\ecE_0)+2\,\deg(\ecE_1)+2\,\deg(\ecE_2)+2\,\deg(\ecE_3)}
\,
q^{\,2\,\deg(\ecE_0)+\sum_i\deg(\ecE_i)-\sum_i \deg(\ecF_i)}
\\
\sum_{\substack{
\bd\in \N^4
\\
~
\\
\psi_1(\bd,\becE)<0
}}
\card{\Homogs_{\frE_0,d_0}}
q^{2(1-g)+2\,d_0+\sumu{i} d_i}
\sum_{(s_i)\in \produ{i} \Homogs_{\frE_i,d_i}}
\!\!\!\!q^{\deg(\pgcd(\ddiv(s_i)+\ecE_i+\ecF_i))}
\,T^{3\,d_0+2\,\sumu{i}d_i},
\end{multline}
est $q^{-1}$-contr\^ol\'ee \`a l'ordre $2$. Le proc\'ed\'e est le m\^eme pour
les s\'eries correspondant aux cinq autres conditions.

D'apr\`es le lemme \ref{lm:ctheta} et la d\'efinition de $\psi_1$,
pour tout r\'eel $\theta>0$,  la s\'erie
$\Zpr{0,\text{err},\psi_1}(\becE ,\bfrE, T)$ 
est major\'ee par
\begin{multline}
c_{\theta}\,c^{1+\sumu{_v} 2\,\Max(v(\ecE_i)+v(\ecF_i))-\Min(\ecE_i+\ecF_i)}
T^{3\,\deg(\ecE_0)+2\,\sumu{i}\deg(\ecE_i)}
\,
q^{2\,\deg(\ecE_0)+\sumu{i}\deg(\ecE_i)-\sumu{i} \deg(\ecF_i)}
\\
\sum_{\substack{
d_1\geq 0
\\
0\leq d_0 < \deg(\ecF_1)
\\
0\leq d_2 < \deg(\ecF_1)
\\
0\leq d_3 < \deg(\ecF_1)
}}
(d_1+1)\,
q^{3+2(1-g)+3\,d_0
+2\,d_1+(2+\theta)\,(d_2+d_3)
}
q^{\deg(\pgcd(\ecE_i+\ecF_i))} 
\,T^{3\,d_0+2\,\sumu{i} d_i},
\end{multline}
donc par 
\begin{equation}
\frac{1}{(1-q\,T)^2}\wtZpr{0,\text{err},\psi_1,\theta}(\becE, T)
\end{equation}
o\`u $\wtZpr{0,\text{err},\psi_1,\theta}(\becE, T)$
est donn\'ee par l'expression
\begin{multline}
c_{\theta}\,c^{1+\sumu{_v} 2\,\Max(v(\ecE_i)+v(\ecF_i))-\Min(v(\ecE_i)+v(\ecF_i))}
\\
\times \,T^{3\,\deg(\ecE_0)+2\,\sumu{i}\deg(\ecE_i)}
\,
q^{2\,\deg(\ecE_0)+\sum_i\deg(\ecE_i)-\sum_i \deg(\ecF_i)}
\\
\times\sum_{\substack{
0\leq d_0 < \deg(\ecF_1)
\\
0\leq d_2 < \deg(\ecF_1)
\\
0\leq d_3 < \deg(\ecF_1)
}}
q^{3+2(1-g)+3\,d_0+(2+\theta)(d_2+d_3)}
q^{\deg(\pgcd(\ecE_i+\ecF_i))} 
\,T^{3\,d_0+2\,d_2+2\,d_3}.
\end{multline}
En raisonnant comme dans la preuve du lemme \ref{lm:convz4}
on montre que pour $\theta>0$ assez petit,
la s\'erie 
\begin{equation}
\sum_{\becE \in
\diveffc^{7}}
\mu_S(\becE)
\wtZpr{0,\text{err},\psi_1,\theta}(\becE, q^{-1+\theta})
\end{equation}
est absolument convergente, ce qui
permet de conclure.
\end{demo}

\subsubsection{Le terme $\Zpr{0,\text{princ}}$}
\begin{prop}\label{prop:z0princ}
La s\'erie 
\begin{equation}
\sum_{\becE \in
\diveffc^{7}}
\mu_S(\becE) 
\sum_{\bfrE\in \piczcSe }
\Zpr{0,\text{princ}}(\becE ,\bfrE, T)
\end{equation}
s'\'ecrit
\begin{equation}
Z_{\courbe}(q^2\,T^3)\,Z_{\courbe}(q\,T^2)^3\,\wtZpr{0,\text{princ}}(T)
\end{equation}
ou $\wtZpr{0,\text{princ}}(T)$ est une s\'erie  de rayon de
convergence
strictement sup\'erieur \`a $q^{-1}$, v\'erifiant
\begin{equation}
\wtZpr{0,\text{princ}}\left(q^{-1}\right)=
(q-1)^4\,q^{\,2(1-g)}
\prod_{v\in\courbe^{(0)}}
(1-q_v^{-1})^{\rg(\Pic(S))}\,\frac{\card{S(\kappa_{v})}}{q_v^{\,\dim(S)}}.
\end{equation}
\end{prop}
\begin{demo}
Rappelons que $\Zpr{0,\text{princ}}(\becE ,\bfrE, T)$
est donn\'ee par l'expression
\begin{multline}
T^{3\,\deg(\ecE_0)+2\,\sumu{i}\deg(\ecE_i)}
\,
q^{\,2\,(1-g)+2\,\deg(\ecE_0)+\sumu{i}\deg(\ecE_i)-\sumu{i}\deg(\ecF_i)}
\\
\times \sum_{
\substack{
\bd\in \N^4
}}
\card{\Homogs_{\frE_0,d_0}}
q^{\,2\,d_0+\sumu{i} d_i}
\sum_{(s_i)\in \produ{i} \Homogs_{\frE_i,d_i}}
\!\!\!\!q^{\deg(\pgcd(\ddiv(s_i)+\ecE_i+\ecF_i))}
\,T^{3\,d_0+2\,\sumu{i}d_i}.
\end{multline}
D'apr\`es la remarque \ref{rq:nis}, pour tout $\bd\in \N^4$, l'application
\begin{equation}
(s_0,(s_i))\mapsto (\ddiv(s_0),(\ddiv(s_i)))
\end{equation}
induit une surjection de l'ensemble
\begin{equation}
\disju{\bfrE\in \piczcSe }
\Homogs_{\bfrE_0,d_0}\times \prod_i \Homogs_{\bfrE_i,d_i}
\end{equation}
sur l'ensemble des \'el\'ements $(\ecG_0,(\ecG_i))\in \diveffc^4$ 
de degr\'e $\bd$, surjection dont les fibres sont de cardinal $(q-1)^4$.
Ainsi la somme
\begin{equation}
\sum_{\bfrE\in \piczcSe }\Zpr{0,\text{princ}}(\becE ,\bfrE, T)
\end{equation}
est \'egale \`a 
\begin{multline}
(q-1)^4\,
T^{3\,\deg(\ecE_0)+2\,\sumu{i}\deg(\ecE_i)}
\,
q^{\,2(1-g)+2\,\deg(\ecE_0)+\sumu{i}\deg(\ecE_i)-\sumu{i}\deg(\ecF_i)}
\\
\times\sum_{
(\ecG_0,(\ecG_i))\in \diveffc^4
}
q^{\,2\,\deg(\ecG_0)+\sumu{i} \deg(\ecG_i)+\deg(\pgcd(\ecG_i+\ecE_i+\ecF_i))}
\,T^{3\,\deg(\ecG_0)+2\,\sumu{i}\deg(\ecG_i)}.
\end{multline}
L'expression pr\'ec\'edente peut s'\'ecrire comme le produit eul\'erien
({\it cf.} les notations \ref{nota:Fnu})
\begin{multline}
(q-1)^4\,q^{\,2(1-g)}\,
\,
T^{3\,\deg(\ecE_0)+2\,\sumu{i}\deg(\ecE_i)}
\,
q^{2\,\deg(\ecE_0)+\sumu{i}\deg(\ecE_i)-\sumu{i}
  \deg(\ecF_i)}
\\
\times
Z_{\courbe}(q^2\,T^3)
\prod_{v\in \courbe^{(0)}}
F_{(v(\ecF_i)+v(\ecE_i))}(q_v,q_v\,T^{2\,f_v})
\end{multline}

On pose, pour tout $v\in \courbe^{(0)}$,
\begin{multline}
Z_{1,v}(T)\eqdef
\sum_{(e_0,\be,\bbf)\in \{0,1\}^7}
\mu^0_S(e_0,\be,\bbf)
\,(q^2\,T^3)^{f_v\,\be_0}\,\,  (q\,T^2)^{f_v\,\sum_{i}e_i}\,\,q_v^{-\sum_{i}f_i}\,\,  
\\
\times 
F_{(e_i+f_i)}(q_v,q_v\,T^{2\,f_v}).
\end{multline}
Ainsi, on a
\begin{align}
\Zpr{0,\text{princ}}(T)&=
\sum_{\becE \in
\diveffc^{7}}
\mu_S(\becE)
\sum_{\bfrE\in \piczcSe }
\Zpr{0,\text{princ}}(\becE ,\bfrE, T)
\\
&=
(q-1)^4\,q^2\, Z_{\courbe}(q^2\,T^3)\prod_{v\in \courbe^{(0)}} Z_{1,v}(T).
\end{align}
On pose
\begin{equation}
Z_{2,v}(T)
\eqdef
\left(1-(q\,T^2)^{f_v}\right)^3
Z_{1,v}(T)
\end{equation}
de sorte qu'on a 
\begin{equation}
\Zpr{0,\text{princ}}=
(q-1)^4 q^{\,2(1-g)}
Z_{\courbe}(q^2\,T^3)\,Z_{\courbe}(q\,T^2)^3
\prod_{v\in \courbe^{(0)}} Z_{2,v}(T).
\end{equation}
En reprenant les notations \ref{nota:Fnu},
$Z_{2,v}(T)$ peut s'\'ecrire
\begin{multline}
\frac{1}{1-(q^4\,T^6)^{f_v}}
\sum_{(e_0,\be,\bbf)\in \{0,1\}^7}
\mu^0_S(e_0,\be,\bbf)
(q^2\,T^3)^{f_v\,e_0}\,\,(q\,T^2)^{f_v\,\sum_{i}e_i} \,\phantom{(}q_v\phantom{)}^{\!\!-\sum_{i}f_i}  
\\
\times\,
\wt{F}_{(e_i+f_i)}(q_v,q_v\,T^{2\,f_v}).
\end{multline}
La proposition \ref{prop:coefFnu} 
montre  que la s\'erie $\produ{v} Z_{2,v}(T)$ a un rayon de convergence
strictement sup\'erieur \`a $q^{-1}$.

Pour terminer la d\'emonstration, il suffit donc de montrer qu'on a pour tout
$v\in \courbe^{(0)}$ la relation 
\begin{equation}\label{eq:z2vq}
Z_{2,v}(q^{-1})=(1-q_v^{-1})^{\rg(\Pic(S))}\,\frac{\card{S(\kappa_{v})}}{q_v^{\,\dim(S)}}.
\end{equation}
On pose, pour $(e_0,\be,\bbf)\in \{0,1\}^7$ 
\begin{equation}
\fact_v(e_0,\be,\bbf)
\eqdef
\frac{1}{1-q_v^{-2}}\,q_v^{-\be_0-\sumu{i} (e_i+f_i)}
\wt{F}_{(e_i+f_i)}\left(q_v,q_v^{-1}\right).
\end{equation}
On a donc 
\begin{equation}
Z_{2,v}(q^{-1})=\sum_{(e_0,\be,\bbf)\in \{0,1\}^7}
\mu^0_S(e_0,\be,\bbf)\,\fact_v(e_0,\be,\bbf).
\end{equation}
Le lemme \ref{lm:rel:mir} ci-dessous
et le lemme \ref{lm:rel:densv}
montrent 
que la relation \eqref{eq:z2vq} est bien v\'erifi\'ee,
ce qui conclut la d\'emonstration.
\end{demo}
\begin{lemme}\label{lm:rel:mir}
On a la relation
\begin{multline}\label{eq:rel:mir}
\sum_{(e_0,\be,\bbf)\in \{0,1\}^7}
\mu^0_S(e_0,\be,\bbf)\,\fact_v(e_0,\be,\bbf)
\\
=
\sum_{(e_0,\be,\bbf)\in \{0,1\}^7}
\mu^0_S(e_0,\be,\bbf)\,\dens_{S,v}(e_0,\be,\bbf).
\end{multline}
\end{lemme}
\begin{demo}
Rappelons qu'on a ({\it cf.} les notations \ref{nota:LA} et \ref{nota:dens})
\begin{equation}\label{eq:densSv}
\dens_{S,v}(e_0,\be,\bbf)=
\frac{
\card{
\left\{(x_0,\bx,\by)\in
  \kappa_v^{(e_0,\be,\bbf)},
\quad \sumu{1\leq i\leq 3}{x_i\,y_i}=0
\right\}
}
}
{q_v^{\,6}}.
\end{equation}

La relation \eqref{eq:rel:mir} (voire m\^eme plus directement le 
fait que le membre de gauche de \eqref{eq:rel:mir}
co\"\i ncide avec le membre de droite de \eqref{eq:z2vq})
peut tr\`es bien se v\'erifier par force brute,
avec l'aide par exemple d'un logiciel de calcul formel.
Montrons comment on peut la retrouver via un minimum de calcul.

Fixons $e_0\in \{0,1\}$.
Tout d'abord, rappelons qu'on a d'apr\`es le lemme \ref{lm:sum:un:imp:mus:zero}
\begin{equation}\label{eq:rel;mu0S}
\sum_{(\be,\bbf)\in \{0,1\}^6}
\mu^0_S(e_0,\be,\bbf)
=0.
\end{equation}
Ensuite, on va montrer ci-dessous qu'on a 
\begin{multline}\label{eq:rel:factv:densv}
\forall (\be,\bbf)\in \{0,1\}^6,\quad
\fact_v(e_0,\be,\bbf)-\fact_v(e_0,1,1)
\\
=\dens_{S,v}(e_0,\be,\bbf)-\dens_{S,v}(e_0,1,1).
\end{multline}
Les deux relations \eqref{eq:rel;mu0S}
et \eqref{eq:rel:factv:densv} montrent le lemme.

Pour montrer la relation \eqref{eq:rel:factv:densv},
on commence par remarquer que par d\'efinition de 
de $\wt{F}_{(e_i+f_i)}$ on a
\begin{align}
\fact_v(e_0,\be,\bbf)
&
=
(1-q_v^{-1})^3\,q_v^{-\be_0-\sum_{i} (e_i+f_i)}
\sum_{\bn\in \N^3}
q_v^{\Min(n_i+e_i+f_i)}\,q_v^{-\sum_{i} n_i}
\\
&
=
(1-q_v^{-1})^3\,
q_v^{\,-e_0}\,
\sum_{\substack{
\bn\in \N^3
\\
n_i\geq e_i+f_i
}}
q_v^{\Min(n_i)}\,q_v^{-\sum_{i} n_i}.
\label{eq:expr:factv}
\end{align}
Par sym\'etrie, il suffit de montrer qu'on a pour tout
$(e_2,e_3,\bbf)\in \{0,1\}^5$ la relation
\begin{multline}\label{eq:rel2}
\fact_v(e_0,(0,e_2,e_3),\bbf)-\fact(e_0,(1,e_2,e_3),\bbf)
\\
=\dens_{S,v}(e_0,(0,e_2,e_3),\bbf)-\dens_{S,v}(e_0,(1,e_2,e_3),\bbf)
\end{multline}
Consid\'erons le cas o\`u $f_1=0$. D'apr\`es \eqref{eq:densSv}, le membre de droite de \eqref{eq:rel2}
vaut alors
\begin{align}
q_v^{-e_0-5}\,\card{\left\{(\bx,\by)\in \kappa_v^{(0,e_2,e_3,\bbf)}, \,\,x_1\neq
  0, \,\,\sum {x_i\,y_i}=0\right\}}
&
=
q_v^{-e_0-5}\,(q_v-1)\card{\kappa_v^{(e_2,e_3,f_2,f_3)}}
\\
&
=
(1-q_v^{-1})\,q_v^{\,-e_0-e_2-e_3-f_2-f_3}
\end{align}
et d'apr\`es \eqref{eq:expr:factv} celui de gauche vaut 
\begin{align}
(1-q_v^{-1})^3\,q_v^{-e_0}
\sum_{\substack{
(n_2,n_3)\in \N^2
\\
n_2\geq e_2+f_2
\\
n_3\geq e_3+f_3
}}
q_v^{\Min(0,n_2,n_3)}\,q_v^{-n_2-n_3}
=
(1-q_v^{-1})\,q_v^{-e_0-e_2-e_3-f_2-f_3} 
\end{align}
d'o\`u l'\'egalit\'e cherch\'ee dans ce cas.

Supposons \`a pr\'esent $f_1=1$. D'apr\`es \eqref{eq:densSv}, le membre de droite de \eqref{eq:rel2} vaut 
\begin{equation}
q_v^{-e_0-5}\,(q_v-1)\card{\left\{(x_2,y_2,x_3,y_3)\in
  \kappa_v^{(e_2,e_3,f_2,f_3)}, 
\,\,x_2\,y_2+x_3\,y_3=0\right\}}.
\end{equation}
Par un calcul facile, on trouve que cette quantit\'e est \'egale \`a
\begin{equation}
\left\{
\begin{array}{ll}
q_v^{\,-1-e_0-e_2-e_3-f_2-f_3}(q_v-1)&\text{si}\quad e_2+f_2\geq
1\quad\text{et}\quad e_3+f_3\geq 1\\
q_v^{\,-3-e_0-e_2-f_2}(q_v-1)(2\,q_v-1)&\text{si}\quad e_2+f_2\geq 1
\quad\text{et}\quad e_3+f_3=0\\
q_v^{\,-5-e_0}\,(q_v-1)(q_v^{3}+q_v^2-q_v)&\text{si}\quad e_2+f_2=0
\quad\text{et}\quad e_3+f_3=0.\\
\end{array}
\right.
\end{equation}
D'apr\`es \eqref{eq:expr:factv} le membre de gauche de \eqref{eq:rel2} vaut 
\begin{equation}\label{eq:mbre:gauche}
(1-q_v^{-1})^3\,q_v^{-e_0}
\sum_{\substack{
(n_2,n_3)\in \N^2
\\
n_2\geq e_2+f_2
\\
n_3\geq e_3+f_3
}}
q_v^{\Min(1,n_2,n_3)}\,q_v^{-1-n_2-n_3}.
\end{equation}
Si $e_2+f_2\geq 1$ et $e_3+f_3\geq 1$, l'expression 
\eqref{eq:mbre:gauche} s'\'ecrit
\begin{equation}
(1-q_v^{-1})^3\,q_v^{-e_0}
\sum_{\substack{
(n_2,n_3)\in \N^2
\\
n_2\geq e_2+f_2
\\
n_3\geq e_3+f_3
}}
\,q_v^{-n_2-n_3}
=
(1-q_v^{-1})q_v^{\,-e_0-e_2-f_2-e_3-f_3}.
\end{equation}
Si $e_2+f_2\geq 1$ et $e_3+f_3=0$, 
on a
\begin{align}
\eqref{eq:mbre:gauche}&=(1-q_v^{-1})^3\,q_v^{-e_0}
\left[
\sum_{\substack{
(n_2,n_3)\in \N^2
\\
n_2\geq e_2+f_2
\\
n_3\geq 1
}}
\,q_v^{-n_2-n_3}
+
\sum_{\substack{
n_2\geq e_2+f_2
}}
\,q_v^{-1-n_2}
\right]
\\
&=
q_v^{-e_0}
\left[
(1-q_v^{-1})\,
q_v^{-e_2-f_2-1}
+(1-q_v^{-1})^2
q_v^{-e_2-f_2-1}
\right]
\\
&
=q_v^{\,-3-e_0-e_2-f_2}
(q_v-1)
(
2\,q_v-1
)
\end{align}
Enfin, si $e_2+f_2=e_3+f_3=0$,
on a 
\begin{align}
\eqref{eq:mbre:gauche}&=(1-q_v^{-1})^3\,q_v^{-e_0}
\left[
q_v^{-1}
+
\sum_{n_2\geq 1}
q_v^{-1-n_2}
+
\sum_{n_3\geq 1}
\,q_v^{-1-n_3}
+
\sum_{\substack{
n_2\geq 1
\\
n_3\geq 1
}}
\,q_v^{\,-n_2-n_3}
\right]
\\
&=
q_v^{-e_0}
\left[
q_v^{-1}(1-q_v^{-1})^3\,
+2\,q_v^{-2}(1-q_v^{-1})^2\,
+q_v^{-2}(1-q_v^{-1})
\right]
\\
&
=q_v^{-5-e_0}
(q_v-1)
\left(
q_v^3+q_v^2-q_v
\right)
\end{align}
d'o\`u le r\'esultat cherch\'e.
\end{demo}

\bibliographystyle{alpha}

\begin{thebibliography}{dlBBD07}

\bibitem[BM90]{BaMa:pts_rat}
V.~V. Batyrev and Yu.~I. Manin.
\newblock Sur le nombre des points rationnels de hauteur born\'e des
  vari\'et\'es alg\'ebriques.
\newblock {\em Math. Ann.}, 286(1-3):27--43, 1990.

\bibitem[Bou03]{Bou:vtetor}
David Bourqui.
\newblock Fonction z\^eta des hauteurs des vari\'et\'es toriques d\'eploy\'ees
  dans le cas fonctionnel.
\newblock {\em J. Reine Angew. Math.}, 562:171--199, 2003.

\bibitem[Bro07]{Br:manin:conjecture:dim:2}
T.D. Browning.
\newblock The {M}anin conjecture in dimension 2.
\newblock Lecture notes for the "School and conference on analytic number
  theory", ICTP, Trieste, 23/04/07-11/05/07 {\url{arXiv:0704.1217v1}}, 2007.

\bibitem[CLT00]{CLT_vect1}
Antoine Chambert-Loir and Yuri Tschinkel.
\newblock Points of bounded height on equivariant compactifications of vector
  groups. {I}.
\newblock {\em Compositio Math.}, 124(1):65--93, 2000.

\bibitem[CLT02]{CLT_vect3}
Antoine Chambert-Loir and Yuri Tschinkel.
\newblock On the distribution of points of bounded height on equivariant
  compactifications of vector groups.
\newblock {\em Invent. Math.}, 148(2):421--452, 2002.

\bibitem[Cox95a]{Cox:funct}
David~A. Cox.
\newblock The functor of a smooth toric variety.
\newblock {\em Tohoku Math. J. (2)}, 47(2):251--262, 1995.

\bibitem[Cox95b]{Cox:hom_coo_ring}
David~A. Cox.
\newblock The homogeneous coordinate ring of a toric variety.
\newblock {\em J. Algebraic Geom.}, 4(1):17--50, 1995.

\bibitem[Der06]{Der:sdp:ut:hyp}
Ulrich Derenthal.
\newblock Singular {D}el {P}ezzo surfaces whose universal torsors are
  hypersurfaces.
\newblock {\url{arXiv:math/0604194v1}}, 2006.

\bibitem[dlB02]{dlB:duke}
R{\'e}gis de~la Bret{\`e}che.
\newblock Nombre de points de hauteur born\'ee sur les surfaces de del {P}ezzo
  de degr\'e 5.
\newblock {\em Duke Math. J.}, 113(3):421--464, 2002.

\bibitem[dlBB07]{dlBBr:sdP4I}
R.~de~la Bret{\`e}che and T.~D. Browning.
\newblock On {M}anin's conjecture for singular del {P}ezzo surfaces of degree
  4. {I}.
\newblock {\em Michigan Math. J.}, 55(1):51--80, 2007.

\bibitem[dlBBD07]{dlBBD}
R{\'e}gis de~la Bret{\`e}che, Tim~D. Browning, and Ulrich Derenthal.
\newblock On {M}anin's conjecture for a certain singular cubic surface.
\newblock {\em Ann. Sci. \'Ecole Norm. Sup. (4)}, 40(1):1--50, 2007.

\bibitem[FMT89]{FMT}
Jens Franke, Yuri~I. Manin, and Yuri Tschinkel.
\newblock Rational points of bounded height on {F}ano varieties.
\newblock {\em Invent. Math.}, 95(2):421--435, 1989.

\bibitem[Has04]{Has:eq:ut:cox:rings}
B.~Hassett.
\newblock Equations of universal torsors and {C}ox rings.
\newblock In {\em Mathematisches {I}nstitut, {G}eorg-{A}ugust-{U}niversit\"at
  {G}\"ottingen: {S}eminars {S}ummer {T}erm 2004}, pages 135--143.
  Universit\"atsdrucke G\"ottingen, G\"ottingen, 2004.

\bibitem[HK00]{hukeel:mori}
Yi~Hu and Sean Keel.
\newblock Mori dream spaces and {GIT}.
\newblock {\em Michigan Math. J.}, 48:331--348, 2000.
\newblock Dedicated to William Fulton on the occasion of his 60th birthday.

\bibitem[Pey03]{Pey:var_drap}
Emmanuel Peyre.
\newblock Points de hauteur born{\'e}e sur les vari{\'e}t{\'e}s de drapeaux en
  caract{\'e}ristique finie.
\newblock {\url{arXiv:math/0303067v1}}, 2003.

\bibitem[Sal98]{Sal:tammes}
Per Salberger.
\newblock Tamagawa measures on universal torsors and points of bounded height
  on {F}ano varieties.
\newblock {\em Ast\'erisque}, (251):91--258, 1998.
\newblock Nombre et r\'epartition de points de hauteur born\'ee (Paris, 1996).

\end{thebibliography}

\end{document}